\documentclass{amsart}

\usepackage{amscd,amstext,amsfonts,amsthm,amsxtra,amsmath}
\usepackage{amssymb}
\usepackage{mathrsfs}
\usepackage[all]{xy}

\usepackage[margin=1in]{geometry}
\theoremstyle{plain}
\newtheorem{theorem}[equation]{Theorem}
\newtheorem{lemma}[equation]{Lemma}
\newtheorem{proposition}[equation]{Proposition}
\newtheorem{corollary}[equation]{Corollary}

\theoremstyle{definition}
\newtheorem{definition}[equation]{Definition}
\newtheorem{example}[equation]{Example}

\theoremstyle{plain}
\newtheorem{remark}[equation]{Remark}

\def\Flag{\operatorname{Flag}}
\def\OG{\operatorname{OG}}
\def\LG{\operatorname{LG}}
\def\Gr{\operatorname{Gr}}
\def\Rad{\operatorname{Rad}}

\def\SL{\operatorname{SL}}
\def\Sp{\operatorname{Sp}}
\def\OO{\operatorname{O}}
\def\SO{\operatorname{SO}}
\def\PO{\operatorname{PO}}
\def\PGL{\operatorname{PGL}}
\def\PSL{\operatorname{PSL}}
\def\PSp{\operatorname{PSp}}

\title[Diagram for VMRTs on the wonderful symmetric varieties]{Diagrams for varieties of minimal rational tangents on the wonderful symmetric varieties}
\author[S.-y. Kim]{Shin-young Kim}
\address{Basic Science Research Institute\\ Ewha Womans University\\ Seoul 03760, Korea}
\email{shinyoung.kim@ewha.ac.kr}
\date{\today}
\begin{document}
\maketitle
\begin{abstract} We describe varieties of minimal rational tangents on the wonderful symmetric varieties by marked Dynkin diagrams. An irreducible component of a variety of minimal rational tangents is a rational homogeneous space, and hence, we have a corresponding marked Dynkin diagram expression. On the other hand, we have marked Dynkin diagrams from the marked Kac diagram of a irreducible symmetric space by marking adjacent nodes to the marked node. We note that the above two diagrams coincide when the restricted root system is not of type A.
\end{abstract}

\section{Introduction}
Let $G$ be a simple Lie group and $T$ be a maximal torus. Let $B$ be a Borel subgroup and $P^{\alpha}$ be a maximal Parabolic subgroup associated with a simple root $\alpha$. Then we get a \emph{marked Dynkin diagram associated with $G/P^{\alpha}$} which is a Dynkin diagram of $G$ with a marking at the simple root $\alpha$. For example, the marked Dynkin diagram associated with a Grassmannian $\Gr(k,n+1)=\PSL_{n+1}\mathbb C /P^{\alpha_k}$ follows.
\setlength{\unitlength}{0.6cm}
    \begin{center}
    \begin{picture}(0,1)(3,0)
        \put(0,0.5){\line(1,0){1}}
        %\put(1,0.5){\line(1,0){.85}}
        \multiput(1,0.5)(0.2,0){5}{\line(1,0){.1}}
        \put(2.15,0.5){\line(1,0){.9}}
        \put(3,0.5){\line(1,0){.9}}
        %\put(4,0.5){\line(1,0){1}}
        \multiput(4,0.5)(0.2,0){5}{\line(1,0){.1}}
        \put(5,0.5){\line(1,0){1}}

        \put(0,0.5){\circle*{.3}}
        \put(1,0.5){\circle*{.3}}
        \put(2,0.5){\circle*{.3}}
        \put(4,0.5){\circle*{.3}}
        \put(5,0.5){\circle*{.3}}
        \put(6,0.5){\circle*{.3}}
        \put(2.7, 0.35){$\times$}
    \end{picture}    
    \end{center}

Assume that $\alpha$ is a long simple root. Then, it is known that (See \cite{HM02}) one can get the marked Dynkin diagram of the variety of minimal rational tangents, variety of lines, of $G/P^{\alpha}$ from the marked Dynkin diagram associated with $G/P^{\alpha}$ by marking the adjacent nodes to the marked node and remove the previous marking and connected edges. For example, the variety of minimal rational tangents $\mathbb P^{k-1} \times \mathbb P^{l-k-1}$ in $\Gr(k,n+1)$ can be expressed by the following marked Dynkin diagram.
\setlength{\unitlength}{0.6cm}
    \begin{center}
    \begin{picture}(0,1)(3,0)
        \put(0,0.5){\line(1,0){1}}
        %\put(1,0.5){\line(1,0){.85}}
        \multiput(1,0.5)(0.2,0){5}{\line(1,0){.1}}
        %\put(2.15,0.5){\line(1,0){.9}}
        %\put(3,0.5){\line(1,0){.9}}
        %\put(4,0.5){\line(1,0){1}}
        \multiput(4,0.5)(0.2,0){5}{\line(1,0){.1}}
        \put(5,0.5){\line(1,0){1}}

        \put(0,0.5){\circle*{.3}}
        \put(1,0.5){\circle*{.3}}
        %\put(2,0.5){\circle*{.3}}
        %\put(4,0.5){\circle*{.3}}
        \put(5,0.5){\circle*{.3}}
        \put(6,0.5){\circle*{.3}}
        
        \put(1.7, 0.35){$\times$}
        \put(3.7, 0.35){$\times$}
    \end{picture}    
    \end{center}

Our question is roughly "Can we describe the variety of minimal rational tangents of a wonderful symmetric variety only from the Kac diagram at a general point?" and we will give an positive answer in this paper, except some cases. When the question is true, our guess is that the irreducible wonderful symmetric variety might be a negative example for the characterization problem: "Is the variety of minimal rational tangents of a wonderful symmetric variety determine the whole variety?". It is because for a symmetric homogeneous space, we need not only the Kac diagram but also the Satake diagram.

Let $X$ be a projective uniruled algebraic variety over complex numbers. An irreducible family of rational curves $\mathcal K$ is called \emph{a covering family} if $\mathcal K_p$ is non-empty for general $p \in X$. A covering family $\mathcal K$ is called \emph{minimal} if $\mathcal K_p$ is projective for general $p \in X$. The image $C$ of parameterized rational curve $f \colon \mathbb P^1 \rightarrow X$ is called \emph{free} if $f^*TX$ is globally generated and the image $C$ is called \emph{embedded} if $f$ is an immersion. The rational map $\tau_p \colon \mathcal K_p \dasharrow \mathbb P(T_p X)$ is defined at an smooth embedded rational curve $C$ which maps to the tangent direction. Assuming that $\mathcal K$ is minimal, then the closure $\mathcal C_p$ of the image of this tangential map $\tau_x$ is called \emph{the variety of minimal rational tangents at $p$}. We refer to \cite{Kollar} and \cite{Hwang} for more details about rational curves and varieties of minimal rational tangents.

Let $G$ be a simple Lie group, $\sigma$ a group involution, and the $H$ a $\sigma$-fixed subgroup. Let $G_{ad}=G/Z(G)$ be the \emph{adjoint group} of $G$. Then, we have isomorphism $G/H = G_{ad}/G^{\sigma}_{ad}$ and we can call $G/H$ \emph{an (irreducible) adjoint symmetric space}. The classification of the irreducible adjoint symmetric space is very classical and goes back to Cartan. If an adjoint symmetric space $G/H$ is given, we get \emph{a marked Kac diagram} (see section 2, \cite[Section 26.5]{Ti11}, and also \cite{Kac}, \cite[Chapter 4.4.7]{OV}). The adjoint symmetric space $G/H$ admits a unique smooth $G$-equivariant embedding $X$ with a simple normal crossing boundary divisor due to \cite{deCP83}. We call $X$ \emph{the wonderful compactification of the symmetric variety $G/H$}. Since $G$ is rational, $X$ is unirational. Let $\mathcal K$ be a family of minimal rational curves on $X$ and $\mathcal C$ a corresponding variety of minimal rational tangents on $X$. By \cite[Theorem 5.1]{BKP}, we see that $\mathcal K_p$ is smooth and the tangential map $\tau_p \colon \mathcal K_p \rightarrow \mathcal C_p$ is an isomorphism at a general point $p \in X$. 

Let $\frak g$ and $\frak h$ be the Lie algebras of $G$ and $H$. We note the decomposition of $\frak g$ into $H$-representations $\frak g= \frak h \oplus \frak p$, where $\frak h=\frak g^{\sigma}$ is the \emph{adjoint representation}, and $\frak p=\frak g^{-\sigma}$ is the \emph{isotropy representation}. Let $x$ be a highest weight vector of the isotropy representation $\mathfrak p$. Then, $Z:=H.[x] \subset \mathbb P \frak p$ is a homogeneous Legendrian submanifold of $M:=G.[x] \subset \mathbb P \frak g$ with respect to the anti-symmetric bilinear form $B_x$ on $\frak g$ defined by $B_x(y,z):=B(x,[y,z])$ (\cite[Proposition 4.34]{BKP}). This anti-symmetric bilinear form is the \emph{Kostant-Kirillov form} after identifying $\frak g$ with the dual $\frak g^{\vee}$. Hence, we call $Z$ \emph{a Legendrian $H$-orbit}.

Then, by \cite[Theorem 1.1]{BF15} and \cite[Proposition 4.34, Theorem 4.42, Theorem 4.43]{BKP}, we have the following:

\begin{theorem}[Theorem \ref{Main}]\label{theorem 1}
Let $X$ be the wonderful compactification of an irreducible adjoint symmetric space $G/H$. Let $\mathcal C_p \subset \mathbb P(T_pX)$ be a variety of minimal rational tangents at the base point $p \in X$. If the restricted root system of $G/H$ is not of type $A$, then $\mathcal C_p \subset \mathbb P(T_pX)$ is isomorphic to a Legendrian H-orbit $Z \subset \mathbb P \frak p$
\end{theorem}

 For the definition of \emph{the restricted root system}, see Section \ref{subsec: diagram for VMRT of type A}. We show that the Legendrian $H$-orbit $Z$ can be recovered from the marked Kac diagram of the irreducible symmetric space $G/H$ in Section 2. Hence, by Theorem \ref{theorem 1} we have the following:

\begin{theorem}[Theorem \ref{thm: main 2}]\label{theorem 2}
      Let $X$ be the wonderful compactification of an irreducible adjoint symmetric homogeneous space $G/H$. Assume that the restricted root system is not of type A. Then, the marked Dynkin diagram of the variety of minimal rational tangents $\mathcal C_p \subset \mathbb P(T_pX)$ at a general point $p$ can be recovered from the marked Kac diagram of the open orbit $G/H$ by marking the adjacent node to a white node as Section \ref{Diagram}.
\end{theorem}

We also state some known facts with diagrams in Section 3 when the restricted root system of $G/H$ is of type $A$. We have a partial result as follows: If the restricted root system of $G/H$ is $A_r$ for $r \geq 2$, then there exists folding from the Dynkin diagram of $G$ to the Dynkin diagram of $H$ and hence, the marked Dynkin diagram of the variety of minimal rational tangents $\mathcal C_p \subset \mathbb P(T_pX)$ is folded into the marked Dynkin diagram of $Z \subset \mathbb P(\mathfrak p)$. Hence, the marked Dynkin diagram of the variety of minimal rational tangents can be recovered from the marked Kac diagram of $G/H$.

\begin{theorem}[Theorem \ref{thm: main 3}]\label{theorem 3}
     Let $X$ be the wonderful compactification of an irreducible adjoint symmetric homogeneous space $G/H$. Assume that the restricted root system of $G/H$ is $A_r$ for $r \geq 2$. Then $G/H$ is either $\PGL_{r+1}\times\PGL_{r+1}/\PGL_{r+1}$, $\SL_{2r+1}/\OO_{2r+1}$, $\SL_{2n}/\Sp_{2n}$ or $E_6/F_4$ and the marked Dynkin diagram of the variety of minimal rational tangents $\mathcal C_p \subset \mathbb P(T_pX)$ at a general point $p$ can be recovered from the marked Kac diagram of the open orbit $G/H$  by marking the adjacent node to a white node as Section \ref{Diagram} and the folding process.
\end{theorem}

\section{Marked Kac diagrams of Symmetric spaces and marked Dynkin diagrams of the H-orbit}

\subsection{Marked Kac diagrams}

Let $G$ be a semi-simple Lie group, $\sigma$ be an involution on $G$. We call the $\sigma$-fixed subgroup $H$ \emph{a symmetric subgroup $H$} of $G$. Let $G_{ad}=G/Z(G)$ be the \emph{adjoint group} of $G$. Then, we have isomorphism $G/H = G_{ad}/G^{\sigma}_{ad}$ and we can call $G/H$ \emph{an adjoint symmetric space}. By section 2.5 of \cite{BKP}, the adjoint symmetric space $G/H$ decomposes into irreducible symmetric spaces of three types; group type, simple type, and Hermitian type. 
\begin{enumerate}
    \item Group type: $(H\times H)/diag(H)$, where $H$ is a simple adjoint group. 

    \item Simple type: $G/H$, where $G$ is simple and the neutral component $H^0$ is simple or decomposes into two simple factors.

    \item Hermitian type: $G/N_G(L)$, where $G$ is simple adjoint and $L\subset G$ is a Levi subgroup.
\end{enumerate}
Moreover, the Hermition type $G/N_G(L)$ is called \emph{exceptional} if $N_G(L)=L$, and $G/N_G(L)$ is called \emph{non-exceptional} otherwise.

From the simple affine Kac-Moody algebra $\frak g$ and the associated affine Dynkin diagram $D_l^{(k)}$ with numbering at each node in \cite[page 44--46]{Kac} where $D_l$ denote the irreducible Dynkin diagram of $\frak g$ and $k=1,2,3$ is the order of an automorphism group of $\frak g$, we can get \emph{the marked Kac diagrams of the irreducible symmetric homogeneous spaces $G/H$} in the following process. If $G/H=(H\times H) / diag (H)$ is of group type, choose the type $D_l$ corresponding $\frak h$ and choose any node of $D_l^{(1)}$ numbered by $1$ and mark the node by white and others by black, then we get the marked Kac diagram of $H\times H / diag (H)$. If $G/H$ is of simple type, there exist two ways. First, choose the type $D_l$ corresponding to the Lie algebra $\frak g$, choose any node of $D_l^{(1)}$ numbered by $2$, and we mark the node by white and others by black. Second, choose the type $D_l$ corresponding to the Lie algebra $\frak g$, choose a node of $D_l^{(2)}$ numbered by $1$, and we mark the node by white and others by black. Then we get the marked Kac diagram of $G/H$. If $G/N_G(L)$ is of Hermitian type, choose the type $D_l$ corresponding to the Lie algebra $\frak g$ and choose two node of $D_l^{(1)}$ numbered by $1$ and mark the nodes by white and others by black, then we get the marked Kac diagram of $G/N_G(L)$. 

Let us define the marked Kac diagram of adjoint symmetric space in another way. Let $T_H$ be a maximal torus of $H$ and let $T=C_G(T_H)$ as the centralizer. Then we have $T^{\sigma, 0}=T_H$ and $T$ is a $\sigma$-stable maximal torus of $G$. We call this $\sigma$-stable maximal torus $T$ \emph{a $\sigma$-stable maximal torus of fixed type}. We can choose a $\sigma$-stable Borel subgroup $B$ of $G$ containing $T$. Let $R$ be the root system of $(G,T)$ and $\triangle$ be the set of simple roots. Then $\sigma$ permutes the simple roots $\triangle$. Let $\alpha$ be a root in $R$ and $e_{\alpha}$ be a root vector in the Lie algebra $\mathfrak g$ of $G$. Let $\bar{\alpha} := \frac{1}{2}(\alpha + \sigma(\alpha))$, $\bar R:=\{ \bar{\alpha} \in \mathfrak X(T) | \alpha \in R \}$ and  $\bar \triangle:=\{ \bar{\alpha} \in \mathfrak X(T) | \alpha \in \triangle \}$. Then it is known that $\bar R$ is a root system with basis $\bar \triangle$. We consider the $H$ isotropy representation $\frak p$. A weight vector of $\frak p$ with respect to the maximal torus $T^{\sigma, 0}$ is one of the forms $2e_{\alpha}$ or $e_{\alpha} - e_{\sigma(\alpha)}$ with the weight $\bar{\alpha}$, depending on whether $\alpha$ is a non-compact imaginary root, that is, $\sigma(\alpha)=\alpha$ and $\sigma(e_{\alpha})=-e_{\alpha}$, or a complex root, that is, $\sigma(\alpha)$ is neither $\alpha$ nor $-\alpha$, respectively. Among the simple roots in $\bar \triangle$, the set of simple roots of $H$ and the lowest weights of $\frak p$ form a affine simple root system with a corresponding affine Dynkin diagram. If we mark the simple roots of $H$ by black and the lowest weights of $\frak p$ by white, we get {\it the marked Kac diagram of $G/H$}. We listed the marked Kac diagrams in Table 1-3 and we follow \cite{bourbaki} the notion of types of $G/H$.

If $G/H$ is of the group type or of the simple type, $\frak p$ is an irreducible representation with respect to $H^0$ and hence with a unique highest weight. So the marked Kac diagram of $G/H$ has a unique white node. When it is of the group type, there exist a unique adjacent node to the white node except $A_n$-type with $n \geq 2$ in which there are two adjacent nodes to the unique white node. When it is of the simple type, there exist a unique adjacent node to the white node except \textrm{BI}, \textrm{CII}, and \textrm{DI} type in which there exist two adjacent nodes.
If $G/H$ is of the Hermitian type, then $H=N_G(L)$ where $L$ is a Levi-subgroup and $\frak p$ decomposes into two irreducible $L$-representations. There exist two highest weight vectors in $\frak p$ up to scalars and hence two closed $L$-orbits in $\mathbb P \frak p$.  We also note that $H^0=L$ and $N_G(L)/L$ are either trivial or a finite group of order 2 generated by $\sigma$. If $G/H$ is exceptional, there are two isomorphic closed $H$-orbits and if $G/H$ is non-exceptional, then the unique $H$-orbit is the union of the two closed $L$-orbits. In particular, we have two white nodes for the Hermitian type. If we choose one white node, there exist the unique adjacent simple root to the white node except AIII type in which there exist two adjacent nodes to the white node.

\begin{table}
\centering
    \begin{tabular}{|c|c|c|} \hline
type & $H$ &  marked Kac diagram\\ \hline
$A_1$&$\PGL_{2}$ & 
\begin{picture}(4,1)(1,0)
\put(2,0.06){\line(1,0){1}}
\put(2,-0.06){\line(1,0){}}

\put(2,0){\circle*{.3}}
\put(3.1,0){\circle{.3}}
\put(2, -.15){$<$}
\put(2.5, -.15){$>$}
\end{picture}
\\
$A_n$&$\PGL_{n+1}$ &  
\begin{picture}(2,2)(2,0)
\put(0,0){\line(1,0){1}}
\put(1,0){\line(1,0){.85}}
\put(2.1,0){\line(1,0){.4}}
\multiput(2.65,0)(0.2,0){5}{\line(1,0){.1}}
\put(3.55,0){\line(1,0){.4}}
\put(4,0){\line(1,0){1}}
\put(5,0){\line(1,0){1}}

\put(0,0){\line(3,1){2.9}}
\put(6,0){\line(-3,1){2.9}}

\put(3,1){\circle{.3}}
\put(0,0){\circle*{.3}}
\put(1,0){\circle*{.3}}
\put(2,0){\circle*{.3}}
\put(4,0){\circle*{.3}}
\put(5,0){\circle*{.3}}
\put(6,0){\circle*{.3}}
\end{picture}

\\

$B_n$ & $\PO_{2n+1}$  & 
\begin{picture}(2,2)(2,0)
\put(0,0){\line(1,0){1}}
\put(1,0){\line(0,1){1}}
\put(1,0){\line(1,0){.85}}
\put(2.1,0){\line(1,0){.4}}
\multiput(2.65,0)(0.2,0){5}{\line(1,0){.1}}
\put(3.55,0){\line(1,0){.4}}
\put(4,0){\line(1,0){1}}
\put(5,0.05){\line(1,0){1}}
\put(5,-0.05){\line(1,0){1}}

\put(0,0){\circle{.3}}
\put(1,0){\circle*{.3}}
\put(1,1){\circle*{.3}}
\put(2,0){\circle*{.3}}
\put(4,0){\circle*{.3}}
\put(5,0){\circle*{.3}}
\put(6,0){\circle*{.3}}

\put(5.2, -.15){$>$}
\end{picture}
\\

$C_n$& $\PSp_{2n}$ & 

\begin{picture}(9,1.5)(-1.5,0)
\put(0,0.05){\line(1,0){1}}
\put(0,-0.05){\line(1,0){1}}
\put(1,0){\line(1,0){.85}}
\put(2.1,0){\line(1,0){.4}}
\multiput(2.65,0)(0.2,0){5}{\line(1,0){.1}}
\put(3.55,0){\line(1,0){.4}}
\put(4,0){\line(1,0){1}}
\put(5,0.05){\line(1,0){1}}
\put(5,-0.05){\line(1,0){1}}

\put(0,0){\circle{.3}}
\put(1,0){\circle*{.3}}
\put(2,0){\circle*{.3}}
\put(4,0){\circle*{.3}}
\put(5,0){\circle*{.3}}
\put(6,0){\circle*{.3}}

\put(0.2, -.15){$>$}
\put(5.2, -.15){$<$}
\end{picture}
\\

$D_n$& $\PO_{2n}$ &  
\begin{picture}(2,1.5)(2,0)
\put(1,0){\line(-1,1){0.7}}
\put(1,0){\line(-1,-1){0.7}}
\put(1,0){\line(1,0){1}}
\put(2.1,0){\line(1,0){.4}}
\multiput(2.65,0)(0.2,0){5}{\line(1,0){.1}}
\put(3.55,0){\line(1,0){.4}}
\put(4,0){\line(1,0){1}}
\put(5,0){\line(1,1){0.7}}
\put(5,0){\line(1,-1){0.7}}

\put(0.3,0.65){\circle{.3}}
\put(0.3,-0.65){\circle*{.3}}
\put(1,0){\circle*{.3}}
\put(2,0){\circle*{.3}}
\put(4,0){\circle*{.3}}
\put(5,0){\circle*{.3}}
\put(5.7,0.65){\circle*{.3}}
\put(5.7,-0.65){\circle*{.3}}

\end{picture}
\\
$E_6$ & $E^{ad}_6$ &  

\begin{picture}(2,3)(3,0)

\put(1,0){\line(1,0){1}}
\put(2,0){\line(1,0){1}}
\put(3,0){\line(1,0){1}}
\put(4,0){\line(1,0){1}}
\put(3,0){\line(0,1){1}}
\put(3,1){\line(0,1){.8}}

\put(1,0){\circle*{.3}}
\put(2,0){\circle*{.3}}
\put(3,0){\circle*{.3}}
\put(4,0){\circle*{.3}}
\put(5,0){\circle*{.3}}
\put(3,1){\circle*{.3}}
\put(3,2){\circle{.3}}
\end{picture}
 \\
$E_7$ & $E^{ad}_7$ &  

\begin{picture}(2,2)(3,0)
\put(1,0){\line(1,0){1}}
\put(2,0){\line(1,0){1}}
\put(3,0){\line(1,0){1}}
\put(4,0){\line(1,0){1}}
\put(5,0){\line(1,0){1}}
\put(6,0){\line(1,0){1}}

\put(4,0){\line(0,1){1}}

\put(1,0){\circle{.3}}
\put(2,0){\circle*{.3}}
\put(3,0){\circle*{.3}}
\put(4,0){\circle*{.3}}
\put(5,0){\circle*{.3}}
\put(6,0){\circle*{.3}}
\put(7,0){\circle*{.3}}
\put(4,1){\circle*{.3}}

\end{picture}

\\
$E_8$ & $E^{ad}_8$ &  

\begin{picture}(10,2)(-2,0)
\put(0,0){\line(1,0){1}}
\put(1,0){\line(1,0){1}}
\put(2,0){\line(1,0){1}}
\put(3,0){\line(1,0){1}}
\put(4,0){\line(1,0){1}}
\put(5,0){\line(1,0){1}}
\put(6,0){\line(1,0){1}}

\put(5,0){\line(0,1){1}}

\put(0,0){\circle{.3}}
\put(1,0){\circle*{.3}}
\put(2,0){\circle*{.3}}
\put(3,0){\circle*{.3}}
\put(4,0){\circle*{.3}}
\put(5,0){\circle*{.3}}
\put(6,0){\circle*{.3}}
\put(7,0){\circle*{.3}}

\put(5,1){\circle*{.3}}

\end{picture}
 \\

$F_4$ & $F^{ad}_4$ & 

\begin{picture}(2,1)(2,0)
\put(0,0){\line(1,0){1}}
\put(1,0){\line(1,0){1}}
\put(3,0){\line(1,0){1}}
\put(2,0.05){\line(1,0){1}}
\put(2,-0.05){\line(1,0){1}}

\put(0,0){\circle{.3}}
\put(1,0){\circle*{.3}}
\put(2,0){\circle*{.3}}
\put(3,0){\circle*{.3}}
\put(4,0){\circle*{.3}}

\put(2.2, -.15){$>$}
\end{picture}
\\
$G_2$ & $G^{ad}_2$ & 

\begin{picture}(4,1)(1,0)
\put(1,0){\line(1,0){1}}
\put(2,0){\line(1,0){1}}
\put(2,0.07){\line(1,0){1}}
\put(2,-0.07){\line(1,0){1}}

\put(1,0){\circle{.3}}
\put(2,0){\circle*{.3}}
\put(3,0){\circle*{.3}}
\put(2.2, -.15){$>$}
\end{picture}   
\\
&& \\
\hline
  \end{tabular}
\caption{Marked Kac diagram of Group type}
\end{table}

\begin{table}
\centering
\begin{tabular}{|c|c|c|c|} 
\hline
type& $G$ & $H$  & marked Kac diagram\\ 
\hline

AI &$\SL_{3}$&$\SO_{3}$  & 
\begin{picture}(4,1)(1,0)
\put(2,0.03){\line(1,0){.9}}
\put(2,-0.03){\line(1,0){.9}}
\put(2,0.1){\line(1,0){.9}}
\put(2,-0.1){\line(1,0){.9}}

\put(2,0){\circle*{.3}}
\put(3,0){\circle{.3}}
\put(2.2, -.15){$<$}
\end{picture} \\

  &$\SL_{2n+1}$&$\SO_{2n+1}$  & 
\begin{picture}(8,1)(-0.5,0.3)
\put(0,0.45){\line(1,0){1}}
\put(0,0.55){\line(1,0){1}}
\multiput(1,0.5)(0.2,0){5}{\line(1,0){.1}}
\put(2,0.5){\line(1,0){1}}
\put(3,0.5){\line(1,0){1}}
\put(4,0.5){\line(1,0){1}}
\multiput(5,0.5)(0.2,0){5}{\line(1,0){.1}}
\put(6,0.45){\line(1,0){1}}
\put(6,0.55){\line(1,0){1}}

\put(0,0.5){\circle*{.3}}
\put(1,0.5){\circle*{.3}}
\put(2,0.5){\circle*{.3}}
\put(3,0.5){\circle*{.3}}
\put(4,0.5){\circle*{.3}}
\put(5,0.5){\circle*{.3}}
\put(6,0.5){\circle*{.3}}
\put(7,0.5){\circle{.3}}

\put(6.3,0.35){$<$}
\put(0.25,0.35){$<$}

\end{picture} 
  
  \\
  &$\SL_{2n}$&$\SO_{2n}$  &
\begin{picture}(8,1.5)(-0.5,0.3)
\put(0.3,1.15){\line(1,-1){0.7}}
\put(0.3,-0.15){\line(1,1){0.7}}
\multiput(1,0.5)(0.2,0){5}{\line(1,0){.1}}
\put(2,0.5){\line(1,0){1}}
\put(3,0.5){\line(1,0){1}}
\put(4,0.5){\line(1,0){1}}
\multiput(5,0.5)(0.2,0){5}{\line(1,0){.1}}
\put(6,0.45){\line(1,0){1}}
\put(6,0.55){\line(1,0){1}}

\put(0.3,1.15){\circle*{.3}}
\put(0.3,-0.15){\circle*{.3}}
\put(1,0.5){\circle*{.3}}
\put(2,0.5){\circle*{.3}}
\put(3,0.5){\circle*{.3}}
\put(4,0.5){\circle*{.3}}
\put(5,0.5){\circle*{.3}}
\put(6,0.5){\circle*{.3}}
\put(7,0.5){\circle{.3}}

\put(6.25,0.35){$<$}
\end{picture} \\

AII &$\SL_{2n}$&$\Sp_{2n}$  & 
\begin{picture}(8,1.5)(-0.5,0.3)
\put(0.3,1.15){\line(1,-1){0.7}}
\put(0.3,-0.15){\line(1,1){0.7}}
\multiput(1,0.5)(0.2,0){5}{\line(1,0){.1}}
\put(2,0.5){\line(1,0){1}}
\put(3,0.5){\line(1,0){1}}
\put(4,0.5){\line(1,0){1}}
\multiput(5,0.5)(0.2,0){5}{\line(1,0){.1}}
\put(6,0.45){\line(1,0){1}}
\put(6,0.55){\line(1,0){1}}

\put(0.3,1.15){\circle{.3}}
\put(0.3,-0.15){\circle*{.3}}
\put(1,0.5){\circle*{.3}}
\put(2,0.5){\circle*{.3}}
\put(3,0.5){\circle*{.3}}
\put(4,0.5){\circle*{.3}}
\put(5,0.5){\circle*{.3}}
\put(6,0.5){\circle*{.3}}
\put(7,0.5){\circle*{.3}}

\put(6.25,0.35){$<$}
\end{picture}\\

\textrm{BI}& $\SO_{2n+1}$ & $\rm S(\OO_m\times \OO_{2n+1-m})$& 

\begin{picture}(8,1.5)(-0.5,0.3)
\put(0.3,1.15){\line(1,-1){0.7}}
\put(0.3,-0.15){\line(1,1){0.7}}
\multiput(1,0.5)(0.2,0){5}{\line(1,0){.1}}
\put(2,0.5){\line(1,0){1}}
\put(3,0.5){\line(1,0){1}}
\put(4,0.5){\line(1,0){1}}
\multiput(5,0.5)(0.2,0){5}{\line(1,0){.1}}
\put(6,0.45){\line(1,0){1}}
\put(6,0.55){\line(1,0){1}}

\put(0.3,1.15){\circle*{.3}}
\put(0.3,-0.15){\circle*{.3}}
\put(1,0.5){\circle*{.3}}
\put(2,0.5){\circle*{.3}}
\put(3,0.5){\circle{.3}}
\put(4,0.5){\circle*{.3}}
\put(5,0.5){\circle*{.3}}
\put(6,0.5){\circle*{.3}}
\put(7,0.5){\circle*{.3}}

\put(6.25,0.35){$>$}
\end{picture}\\

BII & $\SO_{2n+1}$&$\OO_{2n}$ & 

\begin{picture}(8,1.5)(-0.5,0.3)
\put(0.3,1.15){\line(1,-1){0.7}}
\put(0.3,-0.15){\line(1,1){0.7}}
\multiput(1,0.5)(0.2,0){5}{\line(1,0){.1}}
\put(2,0.5){\line(1,0){1}}
\put(3,0.5){\line(1,0){1}}
\put(4,0.5){\line(1,0){1}}
\multiput(5,0.5)(0.2,0){5}{\line(1,0){.1}}
\put(6,0.45){\line(1,0){1}}
\put(6,0.55){\line(1,0){1}}

\put(0.3,1.15){\circle*{.3}}
\put(0.3,-0.15){\circle*{.3}}
\put(1,0.5){\circle*{.3}}
\put(2,0.5){\circle*{.3}}
\put(3,0.5){\circle*{.3}}
\put(4,0.5){\circle*{.3}}
\put(5,0.5){\circle*{.3}}
\put(6,0.5){\circle*{.3}}
\put(7,0.5){\circle{.3}}

\put(6.25,0.35){$>$}
\end{picture}
\\

\textrm{CII}& $\Sp_{2n}$ & $\Sp_{2m}\times \Sp_{2(n-m)}$& 

\begin{picture}(8,1)(-0.5,0.3)
\put(0,0.45){\line(1,0){1}}
\put(0,0.55){\line(1,0){1}}
\multiput(1,0.5)(0.2,0){5}{\line(1,0){.1}}
\put(2,0.5){\line(1,0){1}}
\put(3,0.5){\line(1,0){1}}
\put(4,0.5){\line(1,0){1}}
\multiput(5,0.5)(0.2,0){5}{\line(1,0){.1}}
\put(6,0.45){\line(1,0){1}}
\put(6,0.55){\line(1,0){1}}

\put(0,0.5){\circle*{.3}}
\put(1,0.5){\circle*{.3}}
\put(2,0.5){\circle*{.3}}
\put(3,0.5){\circle{.3}}
\put(4,0.5){\circle*{.3}}
\put(5,0.5){\circle*{.3}}
\put(6,0.5){\circle*{.3}}
\put(7,0.5){\circle*{.3}}

\put(6.3,0.35){$<$}
\put(0.25,0.35){$>$}

\end{picture}
\\

\textrm{DI}& $\SO_{2n}$ &  $\rm S(\OO_{2m+1}\times \OO_{2(n-m)-1})$& 

\begin{picture}(8,1.5)(-0.5,0.3)
\put(0,0.45){\line(1,0){1}}
\put(0,0.55){\line(1,0){1}}
\multiput(1,0.5)(0.2,0){5}{\line(1,0){.1}}
\put(2,0.5){\line(1,0){1}}
\put(3,0.5){\line(1,0){1}}
\put(4,0.5){\line(1,0){1}}
\multiput(5,0.5)(0.2,0){5}{\line(1,0){.1}}
\put(6,0.45){\line(1,0){1}}
\put(6,0.55){\line(1,0){1}}

\put(0,0.5){\circle*{.3}}
\put(1,0.5){\circle*{.3}}
\put(2,0.5){\circle*{.3}}
\put(3,0.5){\circle{.3}}
\put(4,0.5){\circle*{.3}}
\put(5,0.5){\circle*{.3}}
\put(6,0.5){\circle*{.3}}
\put(7,0.5){\circle*{.3}}
\put(6.3,0.35){$>$}
\put(0.25,0.35){$<$}
\end{picture} \\

& &  $\rm S(\OO_{2m}\times \OO_{2(n-m)})$& 
\begin{picture}(8,1.5)(-0.5,0.3)
\put(0.3,1.15){\line(1,-1){0.7}}
\put(0.3,-0.15){\line(1,1){0.7}}
\multiput(1,0.5)(0.2,0){5}{\line(1,0){.1}}
\put(2,0.5){\line(1,0){1}}
\put(3,0.5){\line(1,0){1}}
\put(4,0.5){\line(1,0){1}}
\multiput(5,0.5)(0.2,0){5}{\line(1,0){.1}}
\put(5,0.5){\line(1,0){1}}
\put(6,0.5){\line(1,1){0.7}}
\put(6,0.5){\line(1,-1){0.7}}

\put(0.3,1.15){\circle*{.3}}
\put(0.3,-0.15){\circle*{.3}}
\put(1,0.5){\circle*{.3}}
\put(2,0.5){\circle*{.3}}
\put(3,0.5){\circle{.3}}
\put(4,0.5){\circle*{.3}}
\put(5,0.5){\circle*{.3}}
\put(6,0.5){\circle*{.3}}
\put(6.7,1.15){\circle*{.3}}
\put(6.7,-0.15){\circle*{.3}}

\end{picture}
\\

DII &$\SO_{2n}$&$\OO_{2n-1}$ & 
\begin{picture}(8,1.5)(-0.5,0.3)
\put(0,0.45){\line(1,0){1}}
\put(0,0.55){\line(1,0){1}}
\multiput(1,0.5)(0.2,0){5}{\line(1,0){.1}}
\put(2,0.5){\line(1,0){1}}
\put(3,0.5){\line(1,0){1}}
\put(4,0.5){\line(1,0){1}}
\multiput(5,0.5)(0.2,0){5}{\line(1,0){.1}}
\put(6,0.45){\line(1,0){1}}
\put(6,0.55){\line(1,0){1}}

\put(0,0.5){\circle{.3}}
\put(1,0.5){\circle*{.3}}
\put(2,0.5){\circle*{.3}}
\put(3,0.5){\circle*{.3}}
\put(4,0.5){\circle*{.3}}
\put(5,0.5){\circle*{.3}}
\put(6,0.5){\circle*{.3}}
\put(7,0.5){\circle*{.3}}
\put(6.3,0.35){$>$}
\put(0.25,0.35){$<$}
\end{picture} \\

\textrm{EI}& $E_6$ & $C_4$&

\begin{picture}(2,2)(2.5,-0.3)
\put(1,0){\line(1,0){1}}
\put(2,0){\line(1,0){1}}
\put(3,0.05){\line(1,0){1}}
\put(3,-0.05){\line(1,0){1}}
\put(4,0){\line(1,0){1}}
\put(1,0){\circle*{.3}}
\put(2,0){\circle*{.3}}
\put(3,0){\circle*{.3}}
\put(4,0){\circle*{.3}}
\put(5,0){\circle{.3}}
%\put(3.75,-.15){$\times$}
\put(3.25,-.15){$<$}
\end{picture}
 \\
\textrm{EII}& $E_6$ & $A_5 \times A_1$& 

\begin{picture}(1,3)(3,0)
\put(1,0){\line(1,0){1}}
\put(2,0){\line(1,0){1}}
\put(3,0){\line(1,0){1}}
\put(4,0){\line(1,0){1}}

\put(3,0){\line(0,1){1}}
\put(3,1){\line(0,1){1}}

\put(1,0){\circle*{.3}}
\put(2,0){\circle*{.3}}
\put(3,0){\circle*{.3}}
\put(4,0){\circle*{.3}}
\put(5,0){\circle*{.3}}

\put(3,1){\circle{.3}}
\put(3,2){\circle*{.3}}

\end{picture}
\\

EIV &$E_6$&$F_4$ &  \\

\textrm{EV}& $E_7$ & $A_7$ & 

\begin{picture}(7,2)(0,0.7)
\put(0,1){\line(1,0){1}}
\put(1,1){\line(1,0){.85}}
\put(2.15,1){\line(1,0){.9}}
\put(3,1){\line(1,0){.9}}
\put(4,1){\line(1,0){1}}
\put(5,1){\line(1,0){1}}
\put(3,1){\line(0,1){.9}}

\put(0,1){\circle*{.3}}
\put(1,1){\circle*{.3}}
\put(2,1){\circle*{.3}}
\put(3,1){\circle*{.3}}
\put(4,1){\circle*{.3}}
\put(5,1){\circle*{.3}}
\put(6,1){\circle*{.3}}

\put(3,2){\circle{.3}}

\end{picture}
 \\
\textrm{EVI}& $E_7$ & $D_6 \times A_1$ & 

\begin{picture}(7,2)(0,0.7)
\put(0,1){\line(1,0){1}}
\put(1,1){\line(1,0){.85}}
\put(2.15,1){\line(1,0){.9}}
\put(3,1){\line(1,0){.9}}
\put(4,1){\line(1,0){1}}
\put(5,1){\line(1,0){1}}
\put(3,1){\line(0,1){.9}}

\put(0,1){\circle*{.3}}
\put(1,1){\circle*{.3}}
\put(2,1){\circle*{.3}}
\put(3,1){\circle*{.3}}
\put(4,1){\circle*{.3}}
\put(5,1){\circle{.3}}
\put(6,1){\circle*{.3}}

\put(3,2){\circle*{.3}}

\end{picture}
\\
\textrm{EVIII}& $E_8$ & $D_8$& 

\begin{picture}(7,2)(0,0)
\put(0,0){\line(1,0){1}}
\put(1,0){\line(1,0){1}}
\put(2,0){\line(1,0){1}}
\put(3,0){\line(1,0){1}}
\put(4,0){\line(1,0){1}}
\put(5,0){\line(1,0){1}}
\put(6,0){\line(1,0){1}}
\put(2,0){\line(0,1){0.8}}
\put(0,0){\circle{.3}}
\put(1,0){\circle*{.3}}
\put(2,0){\circle*{.3}}
\put(3,0){\circle*{.3}}
\put(4,0){\circle*{.3}}
\put(5,0){\circle*{.3}}
\put(6,0){\circle*{.3}}
\put(7,0){\circle*{.3}}
\put(2,1){\circle*{.3}}
\end{picture}\\

\textrm{EVIII}& $E_8$ & $E_7\times A_1$& 
\begin{picture}(7,2)(0,0)
\put(0,0){\line(1,0){1}}
\put(1,0){\line(1,0){1}}
\put(2,0){\line(1,0){1}}
\put(3,0){\line(1,0){1}}
\put(4,0){\line(1,0){1}}
\put(5,0){\line(1,0){1}}
\put(6,0){\line(1,0){1}}
\put(2,0){\line(0,1){0.8}}
\put(0,0){\circle*{.3}}
\put(1,0){\circle*{.3}}
\put(2,0){\circle*{.3}}
\put(3,0){\circle*{.3}}
\put(4,0){\circle*{.3}}
\put(5,0){\circle*{.3}}
\put(6,0){\circle{.3}}
\put(7,0){\circle*{.3}}
\put(2,1){\circle*{.3}}
\end{picture}\\

\textrm{FI}& $F_4$ & $C_3\times A_1$& 
\begin{picture}(2,1)(2,-0.3)
\put(1,0){\line(1,0){1}}
\put(2,0){\line(1,0){1}}
\put(3,0.05){\line(1,0){1}}
\put(3,-0.05){\line(1,0){1}}
\put(4,0){\line(1,0){1}}
\put(1,0){\circle*{.3}}
\put(2,0){\circle{.3}}
\put(3,0){\circle*{.3}}
\put(4,0){\circle*{.3}}
\put(5,0){\circle*{.3}}
\put(3.25,-.15){$>$}
\end{picture}\\
\textrm{FII}& $F_4$ & $B_4$& 

\begin{picture}(2,1)(2,-0.3)
\put(1,0){\line(1,0){1}}
\put(2,0){\line(1,0){1}}
\put(3,0.05){\line(1,0){1}}
\put(3,-0.05){\line(1,0){1}}
\put(4,0){\line(1,0){1}}
\put(1,0){\circle*{.3}}
\put(2,0){\circle*{.3}}
\put(3,0){\circle*{.3}}
\put(4,0){\circle*{.3}}
\put(5,0){\circle{.3}}
\put(3.25,-.15){$>$}
\end{picture}\\
\textrm{G}& $G_2$ & $A_1\times A_1$ &

\begin{picture}(4,1)(1,0)
\put(1,0){\line(1,0){0.9}}
\put(2.1,0){\line(1,0){1}}
\put(2.1,0.08){\line(1,0){1}}
\put(2.1,-0.08){\line(1,0){1}}

\put(1,0){\circle*{.3}}
\put(2,0){\circle{.3}}
\put(3,0){\circle*{.3}}
\put(2.2, -.15){$>$}
\end{picture} 

\\&&&\\ \hline
  \end{tabular} \\
\caption {Marked Kac diagram of Simple type}
\end{table}

\begin{table}
\centering
\begin{tabular}{|c|c|c|c|}
\hline
type& $G$ & $L$ &  marked Kac diagram\\ 
\hline

\textrm{AI}&$\PGL_{2}$ &$\PO_{2}$ & 
\begin{picture}(4,1)(1,0)
\put(2,0.06){\line(1,0){.9}}
\put(2,-0.06){\line(1,0){.9}}

\put(2,0){\circle*{.3}}
\put(3,0){\circle{.3}}
\put(2.2, -.15){$<$}
\end{picture}

\\

\textrm{AIII}&$\PGL_{n}$ & $\rm PG(L_{m}\times L_{n-m})$& 
\begin{picture}(8,2)(-1,0)
\put(0,0){\line(1,0){1}}
\put(1,0){\line(1,0){.85}}
\put(2.1,0){\line(1,0){.4}}
\multiput(2.65,0)(0.2,0){5}{\line(1,0){.1}}
\put(3.55,0){\line(1,0){.4}}
\put(4,0){\line(1,0){1}}
\put(5,0){\line(1,0){1}}

\put(0,0){\line(3,1){2.9}}
\put(6,0){\line(-3,1){2.9}}

\put(3,1){\circle{.3}}
\put(0,0){\circle*{.3}}
\put(1,0){\circle*{.3}}
\put(2,0){\circle{.3}}
\put(4,0){\circle*{.3}}
\put(5,0){\circle*{.3}}
\put(6,0){\circle*{.3}}
\end{picture}
\\

\textrm{BI}& $\PO_{2n+1}$&$\rm P(\OO_2 \times \OO_{2n-1})$ &
\begin{picture}(2,2)(2,0)
\put(0,0){\line(1,0){1}}
\put(1,0){\line(0,1){1}}
\put(1,0){\line(1,0){.85}}
\put(2.1,0){\line(1,0){.4}}
\multiput(2.65,0)(0.2,0){5}{\line(1,0){.1}}
\put(3.55,0){\line(1,0){.4}}
\put(4,0){\line(1,0){1}}
\put(5,0.05){\line(1,0){1}}
\put(5,-0.05){\line(1,0){1}}

\put(0,0){\circle{.3}}
\put(1,0){\circle*{.3}}
\put(1,1){\circle{.3}}
\put(2,0){\circle*{.3}}
\put(4,0){\circle*{.3}}
\put(5,0){\circle*{.3}}
\put(6,0){\circle*{.3}}

\put(5.2, -.15){$>$}
\end{picture}

\\

\textrm{CI}& $\PSp_{2n}$&$\PGL_n$ &
\begin{picture}(9,1.5)(-1.5,0)
\put(0,0.05){\line(1,0){1}}
\put(0,-0.05){\line(1,0){1}}
\put(1,0){\line(1,0){.85}}
\put(2.1,0){\line(1,0){.4}}
\multiput(2.65,0)(0.2,0){5}{\line(1,0){.1}}
\put(3.55,0){\line(1,0){.4}}
\put(4,0){\line(1,0){1}}
\put(5,0.05){\line(1,0){1}}
\put(5,-0.05){\line(1,0){1}}

\put(0,0){\circle{.3}}
\put(1,0){\circle*{.3}}
\put(2,0){\circle*{.3}}
\put(4,0){\circle*{.3}}
\put(5,0){\circle*{.3}}
\put(6,0){\circle{.3}}

\put(0.2, -.15){$>$}
\put(5.2, -.15){$<$}
\end{picture}
\\

\textrm{DI}& $\PO_{2n}$&$\rm P(\OO_2 \times \OO_{2(n-1)})$ &
\begin{picture}(2,2)(2,0)
\put(1,0){\line(-1,1){0.7}}
\put(1,0){\line(-1,-1){0.7}}
\put(1,0){\line(1,0){1}}
\put(2.1,0){\line(1,0){.4}}
\multiput(2.65,0)(0.2,0){5}{\line(1,0){.1}}
\put(3.55,0){\line(1,0){.4}}
\put(4,0){\line(1,0){1}}
\put(5,0){\line(1,1){0.7}}
\put(5,0){\line(1,-1){0.7}}

\put(0.3,0.65){\circle{.3}}
\put(0.3,-0.65){\circle{.3}}
\put(1,0){\circle*{.3}}
\put(2,0){\circle*{.3}}
\put(4,0){\circle*{.3}}
\put(5,0){\circle*{.3}}
\put(5.7,0.65){\circle*{.3}}
\put(5.7,-0.65){\circle*{.3}}
\end{picture}
 \\

\textrm{DIII}& $\PO_{2n}$&$\PGL_{n}$ & 
\begin{picture}(2,2)(2,0)
\put(1,0){\line(-1,1){0.7}}
\put(1,0){\line(-1,-1){0.7}}
\put(1,0){\line(1,0){1}}
\put(2.1,0){\line(1,0){.4}}
\multiput(2.65,0)(0.2,0){5}{\line(1,0){.1}}
\put(3.55,0){\line(1,0){.4}}
\put(4,0){\line(1,0){1}}
\put(5,0){\line(1,1){0.7}}
\put(5,0){\line(1,-1){0.7}}

\put(0.3,0.65){\circle{.3}}
\put(0.3,-0.65){\circle*{.3}}
\put(1,0){\circle*{.3}}
\put(2,0){\circle*{.3}}
\put(4,0){\circle*{.3}}
\put(5,0){\circle*{.3}}
\put(5.7,0.65){\circle{.3}}
\put(5.7,-0.65){\circle*{.3}}
\end{picture}

\\

\textrm{EIII}& $E^{ad}_6$ & $D_5 \times \mathbb C^*$ & 
\begin{picture}(2,3)(3,0)

\put(1,0){\line(1,0){1}}
\put(2,0){\line(1,0){1}}
\put(3,0){\line(1,0){1}}
\put(4,0){\line(1,0){1}}
\put(3,0){\line(0,1){1}}
\put(3,1){\line(0,1){.8}}

\put(1,0){\circle{.3}}
\put(2,0){\circle*{.3}}
\put(3,0){\circle*{.3}}
\put(4,0){\circle*{.3}}
\put(5,0){\circle*{.3}}
\put(3,1){\circle*{.3}}
\put(3,2){\circle{.3}}
\end{picture}
 \\

\textrm{EVII}& $E^{ad}_7$ & $E_6 \times \mathbb C^*$ & 
\begin{picture}(2,2)(3,0)
\put(1,0){\line(1,0){1}}
\put(2,0){\line(1,0){1}}
\put(3,0){\line(1,0){1}}
\put(4,0){\line(1,0){1}}
\put(5,0){\line(1,0){1}}
\put(6,0){\line(1,0){1}}

\put(4,0){\line(0,1){1}}

\put(1,0){\circle{.3}}
\put(2,0){\circle*{.3}}
\put(3,0){\circle*{.3}}
\put(4,0){\circle*{.3}}
\put(5,0){\circle*{.3}}
\put(6,0){\circle*{.3}}
\put(7,0){\circle*{.3}}
\put(4,1){\circle{.3}}
\end{picture}
\\
&&& \\
\hline
    \end{tabular}
    \caption {Marked Kac diagram of Hermitian type}
\end{table}

\subsection{The $H$-orbit $Z$}

Let $G/H$ be an irreducible adjoint symmetric space. We consider the marked Kac diagram of $G/H$ and the $H$-orbit $Z=H.[x] \subset \mathbb P \frak p$ at a highest weight vector $x \in \frak p$. In particular, we want to describe the homogeneous space $Z=H/H_{[x]}$ where $H_{[x]}$ is the isotropy subgroup of $H$ at $[x] \in \mathbb P \frak p$.

When $G/H$ is of the group type, so that $\mathfrak{p}$ is the adjoint representation of $H$ and $x$ is a highest root vector, then the next Lemma \ref{isotoropy group type} follows essentially from \cite[Lemma 3.3]{BF15}. We give a complete proof and statement because the proof is more suitable to describe $Z$. We also use the arguments in Lemma \ref{isotoropy simple type} and Lemma \ref{isotoropy hermitian type}.

\begin{lemma}\label{isotoropy group type} Let $G/H$ be an irreducible adjoint symmetric space of the group type except when $H=\PGL_2(\mathbb C)$. Let $x$ be a highest weight vector of $H$-representation $\frak p$. Then there exists a unique white node in the Kac diagram of $G/H$ and the isotropy subgroup of $H$ at $[x] \in \mathbb P \frak p$ is isomorphic to (1) the maximal parabolic subgroup of $H$ associated to the unique adjacent simple root $\alpha$ to the white node except in $A$-type. When $G/H$ is of the $A$-type, the isotropy subgroup of $H$ at $[x]$ is isomorphic to (2) the intersection of two maximal parabolic subgroups associated with each of two adjacent simple roots to the white node.
\end{lemma}

\begin{proof}
We note that $H$ is a simple group. By Borel's fixed point theorem, the isotropy subgroup at the highest weight vector $[x] \in \mathbb P \frak p$ is a parabolic subgroup $P$ of $H$. A conjugation of $H$ that sends the Borel subgroup of $H$ to the opposite Borel subgroup induces a map on $\frak h$ which sends the highest weight vector $x$ to the lowest weight vector $y$. The opposite parabolic subgroup $Q$ of $P$ satisfies $H.[y]=H/Q$. Since the lowest weight of $y$ under $H$-representation on $\frak p$ is given by a white node in Kac diagram, let $\beta$ be the simple root of $G$ corresponding to the white node.

There exists one adjacent node to the white node except in A-type. For the A-type, there are two adjacent nodes to the white node. We call the corresponding simple root(s) of $H$ as the adjacent root(s). Let $\triangle_1=\{\alpha\}$ be the set of adjacent simple roots. If we choose an element $E \in \mathfrak t$ such that $\alpha(E)=1$ for $\alpha \in \triangle_1$ and $\beta(E)=0$ for $\beta \in \triangle_1^c$ and set $R_k=\{ \alpha \in R \mid \alpha(E) =k\}$ for $k \in \mathbb Z$, then $\frak{h}_0 = \frak t \oplus \bigoplus_{\alpha \in R_0} \frak h_{\alpha}$ and
$\frak{h}_k = \bigoplus_{\alpha \in R_k} \frak h_{\alpha}$ for $k\neq 0$. This give us a $\mathbb Z$-gradation on $\frak h = \bigoplus_{i\in \mathbb Z} \frak h_i$ associated with the set of adjacent simple roots $\triangle_1=\{\alpha\}$. To complete the proof, it suffices to show that the Lie algebra $\frak{h}_{[y]}$ of $Q$ is isomorphic to $\bigoplus_{i\leq 0} \frak h_i$.

First, we will check $\bigoplus_{i\leq 0} \frak h_i \subset \frak{h}_{[y]}$. For any negative root $\gamma$ of $H$, we have $[e_\gamma, y]=0$. Since  $\sigma$ permutes the negative roots, $\frak g_{\gamma}$ is in the isotropy Lie algebra $\frak{h}_{[y]}$. Let $\frak t$ be the Lie algebra of the maximal torus $T^{\sigma, 0}$ of $H$. For $t \in \frak t$, $[t, y]=0$. Hence, $\frak t$ is in the isotropy Lie algebra $\frak{h}_{[y]}$. Let $\{ \alpha_i \}$, $i \in I$, be the set of simple roots of $H$ and let $I_H \subset I$ be a sub-index set such that that $\alpha_i$, $i \in I_H$, is a nonzero simple root of $H$ not in the adjacent set $\{\alpha\}$. Let $\alpha_i$, $i \in I_H$. If we assume that $\bar \gamma_i=\bar \beta +\bar \alpha_i$ is a weight of the representation with the corresponding root $\gamma_i=\beta +\alpha_i$, then, since $\langle \alpha_i | \beta \rangle=0$, we have $\langle \gamma_i | \beta \rangle=\langle \beta |\beta \rangle=2$ which is a contradiction. Hence, $[e_{\alpha_i}, e_\beta]=0$ and $[e_{\alpha_i},d\sigma(e_\beta)]=0$. Hence, $\frak g_{\alpha_i}$ is in the isotropy Lie algebra. Similarly, for a positive root $\lambda=\sum_{i \in I_H} c_i \alpha_i$ of $H$ with $c_i \geq 0$, we have $\frak g_{\lambda} \subset \frak h_{[y]}$. Hence, $\bigoplus_{i\leq 0} \frak h_i \subset \frak{h}_{[y]}$.

For a adjacent root $\alpha$ of $H$, we have $\langle \bar\beta|\bar\alpha \rangle= \langle \beta | \alpha \rangle < 0$. Hence, $\gamma= \beta + \alpha$ is a root of $G$ and $[e_{\alpha}, e_\beta-d\sigma(e_\beta)]=e_{\gamma}-d\sigma(e_{\gamma})$ which implies that $\bar\gamma=\bar \beta +\bar \alpha$ is a weight of $\frak p$. Hence $\frak g_{\alpha}$ is not in the Lie algebra $\frak{h}_{[y]}$. Moreover, $\bigoplus_{i > 0} \frak h_i$ is not in $\frak{h}_{[y]}$. In conclusion, $\frak{h}_{[y]}=\bigoplus_{i\leq 0} \frak h_i$.
\end{proof}

\begin{remark}\label{contact} If $\mathfrak{p}$ is the adjoint representation $\mathfrak{h}$ of $H$ and $x$ is a highest root vector, then $Z$ is a standard contact manifold and associated Kac diagram is also known as extended Dynkin diagram. See \cite[Chapter 4.2.5]{OV} or \cite[Section 4.2 and 4.3]{Ya93}. Moreover, the gradation on the Lie algebra $\frak h$ in the proof of Lemma \ref{isotoropy group type} is the contact gradation, that is, the depth $\mu$ is 2, so
$$\frak h = \frak h_{-2}\oplus \frak h_{-1}\oplus \frak h_{0}\oplus \frak h_{1}\oplus \frak h_{2},$$
and $\dim \frak h_{-2}=\dim \frak h_2 =1$. If in addition $H=\PGL_{n+1} \mathbb C$ for $n>2$, then $Z$ is the partial flag manifold $\Flag(1,n-1)$.
\end{remark}

\begin{lemma}\label{isotoropy simple type} Let $G/H$ be an irreducible adjoint symmetric space of the simple type. Let $x$ be a highest weight vector in $\frak p$ of $H^0$-representation. Then there exist a unique white node in the Kac diagram of $G/H$ and the isotropy subgroup of $H^0$ at $[x] \in \mathbb P \frak p$ is isomorphic to either (1) the maximal parabolic subgroup associated to the unique adjacent simple root $\alpha$ to the white node or (2) the product of two maximal parabolic subgroups associated with each of two adjacent simple roots to the white node.
\end{lemma}

\begin{proof}
We note that $H^0$ is simple or decomposes into two simple factors $H^0_1$ and $H^0_2$. When $H^0$ is simple, there is a unique adjacent root $\alpha \in R(H^0)$, where $R(H^0)$ is the root system of $H^0$. When $H^0=(H^0_1 \times H^0_2)/\Gamma$ the group quotient by a finite group $\Gamma$, there are two adjacent roots $\alpha_1 \in R(H^0_1)$ and $\alpha_2 \in R(H^0_2)$, where $R(H^0_1)$ and $R(H^0_2)$ are root systems of $H^0_1$ and $H^0_1$ respectively.

For the first case, the same argument of Lemma \ref{isotoropy group type} works for $H^0$ with one simple root $\alpha$.

For the second case, by Borel's fixed theorem and the same first argument in the proof of Lemma \ref{isotoropy group type}, the isotropy subgroup at the highest weight vector $[x] \in \mathbb P (\frak p)$ is a parabolic subgroup $P$ of $H^0$ with the opposite parabolic subgroup $Q$ and there exists a lowest vector $y$ that satisfies $H.[y]=H/Q$. On the other hand, under the action of $H^0_i$, for $i=1,2$, the isotropy subgroup at the highest weight vector $[x]$ is a parabolic subgroup $P_i$ with the opposite parabolic subgroup $Q_i$. By the same arguments as Lemma \ref{isotoropy group type} for $H^0_i$ with each simple root $\alpha_i$, $P_i$ is the maximal parabolic subgroup of $H^0_i$ associated with the adjacent simple root $\alpha_i \in R(H^0_i)$, for $i=1,2$. Hence, $P=P_1 P_2 \subset H^0$.
\end{proof}

\begin{lemma}\label{isotoropy hermitian type} Let $G/H$ be an irreducible adjoint symmetric space of the Hermitian type. Let $L$ be the Levi subgroup such that $H=N_G(L)$ and $L'$ be the commutator subgroup of $L$. Let $x$ be a highest weight vector in $\frak p$ of $L'$-representation. Then there is a white node in the Kac diagram of $G/H$ such that the isotropy subgroup of $L'$ at $[x] \in \mathbb P \frak p$ is isomorphic to (1) the maximal parabolic subgroup of $L'$ associated to the unique adjacent simple root $\alpha$ to the white node except AIII type. When $G/H$ is of the AIII type, the isotropy subgroup of $H$ at $[x]$ is isomorphic to (2) the product of two maximal parabolic subgroups associated with each of two adjacent simple roots to the white node.
\end{lemma}

\begin{proof}
We note that $L'$ is simple or decomposes into two simple factors. Since $\frak p$ decomposes into two irreducible $L'$-representations, for an irreducible $L'$-representation we can apply the same argument as Lemma \ref{isotoropy simple type} to get the conclusion.
\end{proof}

\subsection{Marked Dynkin diagrams of $Z$} \label{Diagram}

We continue to use the definitions and notations in the previous subsection; an irreducible adjoint symmetric space $G/H$, $H$-representation $\frak p$ and a closed $H$-orbit $Z \subset \mathbb P \frak p$.
The neutral component $H^0$ of $H$ decomposes into the solvable radical $\Rad H^0$ and the commutator subgroup ${H^0}'$. We note that $H^0$ is reductive for all three types and ${H^0}'$ is either simple or decomposes into two simple factors.
Since $Rad(H^0)$ is a central subgroup of $H^0$ which acts $x \in \frak p$ by a scalar, it acts trivially on $H^0.[x] \in \mathbb P \frak p$. Hence, it is enough to consider ${H^0}'$ action at $[x] \in \mathbb P \frak p$ to describe $Z=H/H_{[x]}$. Hence, by Lemma \ref{isotoropy group type}, Lemma \ref{isotoropy simple type}, Lemma \ref{isotoropy hermitian type} and their proofs, we get a diagram expression of $Z=H/H_{[x]}$ from the marked Kac diagram of the symmetric space $G/H$.

\begin{definition} By the following process from (a) to (d), we get marked Dynkin diagrams from the marked Kac diagram of the symmetric space $G/H$.
 \begin{itemize}
     \item[(a)] We choose the white node corresponding to the lowest weight $\beta$ such that $\lambda=w_0 \beta$ is a highest weight with a highest weight vector $x$, where $w_0$ is the unique longest element of the Weyl group which maps $R^+$ to $R^-$.
     
     \item[(b)] After marking the adjacent nodes to the chosen white node, remove all white nodes and remove all edges connected to the white node. If the number of edges from the chosen white node(long weight) to the newly marked adjacent node(short weight) is $n \geq 2$, put numbering $n$ on the marked node of a obtained marked Dynkin diagram. 
     
     \item[(c)] In addition, if $G/H$ is of the Hermitian type, then we choose the other white node (equivalently, the other highest weight) and do the same process as above to get the second marked Dynkin diagram.
     
     \item[(d)] If $G/H$ is of the Hermitian non-exceptional type, then put a double-head arrow between the two marked Dynkin diagrams to express $\sigma$.
 \end{itemize}
For the highest weight vector $x$ and the homogeneous manifold $Z=H.[x] \subset \mathbb P \frak p$, we define {\it the marked Dynkin diagram of $Z$} as follows; if $G/H$ is of the group type or of the simple type, from (a) and (b), we get a marked Dynkin diagram corresponding to $Z$. If $G/H$ is of the Hermitian exceptional type, we get two marked Dynkin diagrams from the above process. Each one corresponds to the $H$-orbit of the corresponding highest weight. If $G/H$ is of the Hermitian non-exceptional type, we get two marked Dynkin diagrams with an arrow which together correspond to $Z$.\end{definition}

\begin{example} The marked Kac diagram of symmetric space $\SL_3/\SO_3$ is 
\setlength{\unitlength}{0.6cm}
\begin{picture}(1.5,1)(0,0.8)
\put(0,1.15){\line(1,0){1}}
\put(0,1.05){\line(1,0){1}}
\put(0.25,0.85){$\bf{>}$}
\put(0,0.95){\line(1,0){1}}
\put(0,0.85){\line(1,0){1}}
\put(1.1,1){\circle{.3}}
%\put(0.8,0.8){$\times$}
\put(0,1){\circle*{.3}}
\end{picture}. Hence, the marked Dynkin diagram of $Z$ is
\begin{picture}(0.5,0.5)(0,0.2)
\put(0.1,0.4){\small{4}}
\put(0,0){$\times$}
\end{picture} corresponding to $\nu_4(\mathbb P^1)$. The marked Kac diagram of symmetric space $\SL_{2n+1}/\SO_{2n+1}$ is
\begin{picture}(5.5,0)(-0.2,-0.1)
\put(0,0.05){\line(1,0){1}}
\put(0.25,-0.15){$\bf{>}$}
\put(0,-0.05){\line(1,0){1}}
\put(1,0){\line(1,0){1}}
\multiput(2,0)(0.2,0){5}{\line(1,0){.1}}
\put(3,0){\line(1,0){1}}
\put(4,0.05){\line(1,0){1}}
\put(4.25,-0.15){$\bf{>}$}
\put(4,-0.05){\line(1,0){1}}
\put(-0.1,0){\circle{.3}}
%\put(-0.2,-0.2){$\times$}
\put(1,0){\circle*{.3}}
\put(2,0){\circle*{.3}}
\put(3,0){\circle*{.3}}
\put(4,0){\circle*{.3}}
\put(5,0){\circle*{.3}}
\end{picture}, and hence, the marked diagram of $Z$ is
\begin{picture}(4.5,0)(0.8,0.2)
\put(1,0){\line(1,0){1}}
\multiput(2,0)(0.2,0){5}{\line(1,0){.1}}
\put(3,0){\line(1,0){1}}
\put(4,0.05){\line(1,0){1}}
\put(4.25,-0.15){$\bf{>}$}
\put(4,-0.05){\line(1,0){1}}
\put(0.9,0.3){\small{2}}
\put(0.8,-0.2){$\times$}
\put(2,0){\circle*{.3}}
\put(3,0){\circle*{.3}}
\put(4,0){\circle*{.3}}
\put(5,0){\circle*{.3}}
\end{picture}, which represent $\nu_2(Q_{2n-1})$.
\end{example}

\section{Marked Dynkin diagrams of VMRTs}
An irreducible family of rational curves $\mathcal K$ is called \emph{a covering family} if $\mathcal K_p$ is non-empty for general $p \in X$. A covering family $\mathcal K$ is called \emph{minimal} if $\mathcal K_p$ is projective for general $p \in X$.
The image $C$ of parameterized rational curve $f \colon \mathbb P^1 \rightarrow X$ is called \emph{free} if $f^*TX$ is globally generated which is a smooth point in the covering family $\mathcal K$. Moreover, the image $C$ is called \emph{embedded} if $f$ is an immersion. The rational map $\tau_p \colon \mathcal K_p \dasharrow \mathbb P(T_p X)$ is defined at an smooth embedded rational curve $C$ which maps to the tangent direction. Assuming that $\mathcal K$ is minimal, then the closure $\mathcal C_p$ of the image of this tangential map $\tau_x$ is called \emph{the variety of minimal rational tangents at $p$}. We refer to \cite{Kollar} and \cite{Hwang} for more details about rational curves and varieties of minimal rational tangents.

Now, we assume that $G$ is a simple adjoint group and $X$ is the wonderful compactification of an irreducible adjoint symmetric homogeneous space $G/H$. In this subsection, we will define {\it the marked Dynkin diagram of variety of minimal rational tangents $\mathcal C_p \subset \mathbb P(T_pX)$} at a general $p \in X$ as Definition \ref{vmrt Legendrian} and Definition \ref{vmrt L of A}.

\subsection{Marked Dynkin diagrams of VMRTs whose restricted root system is not of type A}\label{subsec: diagram for VMRT of type A}
Let $T$ be a $\sigma$-stable torus such that $T^{-\sigma,0}=\{ t \in T | \sigma(t)=t^{-1} \}^0$ is maximal and let $R$ be a root system with base $\triangle$. We call $T$ is \emph{a maximal torus of split type}. Let $\Theta_s$ be the maximal positive root corresponding to the root system $(R,T, \triangle)$ and $\Theta_s^{\vee}$ be the coroot. Let $\Tilde{\alpha}=\alpha-\sigma(\alpha)$ for $\alpha \in R $ and define 
$$\Tilde{R}=\{\Tilde{\alpha} \mid \alpha \in R \}.$$
Then $\Tilde{R}$ is a root system which we call it \emph{the restricted root system}.

Let $P$ be the standard parabolic subgroup of $G$ such that $G/P$ is the unique closed $G$-orbit of the wonderful compactification $X$. Let $P=LP^u$ be a Levi decomposition where $L$ is the Levi factor and $P^u$ is the unipotent radical. Let $\rho_G$ be the sum of all fundamental weights of the root system for $G$, $\rho_L$ be the sum of all fundamental weights of the root system for the Levi factor $L$ and let $\rho_P$ be the half sum of roots of $P^u$. Then we have $ \rho_G=  \rho_L +  \rho_P$ and we define the weight $$\kappa=\sum_{\alpha \in R^+(L)} \alpha=2 \rho_G-2 \rho_P.$$ The $H$-orbit is $Z=H.[x] \subset \mathbb P \frak p$ at a highest weight vector $x \in \frak p$ as Section 2.1. Let $T_H=T^{\sigma,0}$ be a torus of $H$ and $B_H$ a Borel subgroup of $H$. An irreducible curve $C$ through a general $p$ which is stable by $B_H$ is called \emph{a highest weight curve}. The highest weight curve $C$ is a parameterized embedded free rational curve $f:\mathbb P^1 \rightarrow X$ with image $C \subset X$ through $p$ of direction $df(p)=x$ with $T_H$-weight $\lambda$.

\begin{proposition}\cite[Corollary 2.14]{BKP}
 Let $C$ a highest weight curve of weight $\lambda$, then there exist a root $\alpha$ such that $\lambda = \alpha|_{T_H}$ and 
    $C = \overline{U_{\alpha}\cdot x}$.
\end{proposition}

\begin{proposition}\cite[Theorem 4.26]{BKP}\label{prop:dim of vmrt} Let $p\in X$ be a general point and let $C$ be a highest weight curve such that $p \in C$. Then
\begin{eqnarray*}
    \dim \mathcal C_p = \frac{1}{2}<\kappa, \Theta^{\vee}_s-\sigma(\Theta^{\vee}_s)> + \partial X \cdot C  -2 ,
\end{eqnarray*}
and
 $$ <\kappa, \Theta^{\vee}_s-\sigma(\Theta^{\vee}_s)> = 2 (\dim H.[x] + 1 ).$$
\end{proposition}

\begin{proposition}\cite[Remark 4.27]{BKP}\label{1} Let $C$ be a highest weight curve and $p \in C$. Then $\partial X \cdot C=1$ except when the restricted root system is of type $A$ in which $\partial X \cdot C=2$.
\end{proposition}
\begin{proof}
See \cite[Remark 4.27]{BKP} and also, \cite[Proposition 3.8]{BKP}  for group type, \cite[Proposition 3.10]{BKP} for Hermitian type, and \cite[Proposition 3.14]{BKP} for simple type.
\end{proof}

\begin{corollary} The variety of minimal rational tangents $\mathcal C_p$ equals $Z=H.[x] \subset \mathbb P \frak p$ except when the restricted root system is of type $A$.

\end{corollary}
\begin{proof}  The highest weight vector $x$ is the tangent direction of the highest curve $C$ at $p$. Since $H$ is in an isotropy subgroup at $p$, $H.[x] \subset \mathcal C_p$. We have $\partial X \cdot C=1$ by Proposition \ref{1}. Hence, by Proposition \ref{prop:dim of vmrt}, it follows
\begin{eqnarray*}
    \dim \mathcal C_p &=& \frac{1}{2}<\kappa, \Theta^{\vee}_s-\sigma(\Theta^{\vee}_s)> + \partial X \cdot C  -2 \\
     &=&\dim H.[x]+\partial X \cdot C  -1 =\dim H.[x]. \end{eqnarray*}
Hence, $\mathcal C_p$ is $Z=H.[x] \subset \mathbb P \frak p$.
\end{proof}

By \cite[Proposition 4.23]{BKP}, the orbit $M=G.[x] \subset \mathbb P \frak g$ is either the minimal orbit $\mathcal O_{min}(G)$ or the nilpotent orbit $\mathcal O=G.[y]$ where $y=x_{\alpha_1}+x_{\alpha_2} \in \frak g$ is a sum of two non-zero root vectors such that $\alpha_1$ and $\alpha_2$ are strongly orthogonal long roots. If $\partial X \cdot C=2$, from the dimension formula $\dim \mathcal C_p =\dim H.[x]+1$, the next candidate for $\mathcal C_p$ is either $\mathcal O_{min}(G)$ or $\mathcal O$.

To prove Theorem \ref{Main}, we need \cite[Theorem 1.1]{BF15} and \cite[Proposition 4.34, Theorem 4.42, Theorem 4.43]{BKP}. It follows.

\begin{theorem}[Theorem \ref{theorem 1}]\label{Main}
Let $X$ be the wonderful compactification of an irreducible adjoint symmetric space $G/H$. Let $\mathcal C_p \subset \mathbb P(T_pX)$ be a variety of minimal rational tangents at a general point $p \in X$.
\begin{itemize}
    \item[(1)] If $G/H$ is of the group type but not of type $A$,
    then $\mathcal C_p\subset \mathbb P(T_pX)$ is isomorphic to the adjoint variety $X_{ad}^H \subset \mathbb P \frak h$, the Legendrian H-orbit.
    \item[(2)] If $G/H$ is of the simple type or the Hermitian type, whose restricted root system is not of type $A$, then $\mathcal C_p \subset \mathbb P(T_pX)$ is isomorphic to a Legendrian H-orbit $Z \subset \mathbb P \frak p$.
    \item[(3)] If the restricted root system of $G/H$ is of type $A_r$ with $r\geq 2$, then $\mathcal C_p \subset \mathbb P(T_pX)$ is isomorphic to a rational homogeneous space $G/P_{\lambda} \subset \mathbb P V_{\lambda}$ with the highest weight $\lambda$ listed in section \ref{A type}.
\end{itemize}
 Moreover, there exists a unique variety of minimal rational tangents except in the Hermitian exceptional type. The unique variety of minimal tangents has one component for group type and simple type, and two isomorphic components for the Hermitian non-exceptional type.
\end{theorem}

\begin{proof}
When $G/H$ is of the group type whose restricted root system is not of type $A$, by \cite[Theorem 1.1]{BF15}, there is the unique variety of minimal rational tangents $\mathcal C_p\subset \mathbb P(T_pX)$ which is isomorphic to the adjoint variety $X_{ad}^H \subset \mathbb P \frak h$ which is the Legendrian $H$-orbit $Z \subset \mathbb P \frak p$. 
By \cite[Theorem 4.42]{BKP}, the tangent map $\mathcal K_p \rightarrow \mathcal C_p$ is an isomorphism. Hence, by \cite[Proposition 4.34, Theorem 4.43 (3)]{BKP}, there is a $G$-invariant contact structure on a nilpotent $G$-orbit $M \subset \mathbb P \frak g$ which contains a variety of minimal rational tangents $\mathcal C_p \subset \mathbb P(T_pX)$ as a homogeneous Legendrian submanifold isomorphic to $Z \subset \mathbb P \frak p$ when $G/H$ is of the simple type or the Hermitian type whose restricted root system is not of type $A$. If the restricted root system of $G/H$ is of type $A_r$ with $r\geq 2$, by \cite[Theorem 4.43(4)]{BKP} (see also Lemma \ref{lem: A_r}), $\mathcal C_p \subset \mathbb P(T_pX)$ is isomorphic to a rational homogeneous space $G/P_{\lambda} \subset \mathbb P V_{\lambda}$ with the highest weight $\lambda$ listed in section \ref{A type}. The last statement also follows by \cite[Theorem 4.43(3)]{BKP}.
\end{proof}

\begin{definition}\label{vmrt Legendrian}
Let $X$ be the wonderful compactification of an irreducible adjoint symmetric homogeneous space $G/H$ whose restricted root system is not of type $A$. Let $\mathcal C_p \subset \mathbb P(T_pX)$ be a variety of minimal rational tangents at a general point $p \in X$ which corresponds to $Z \subset \mathbb P \frak p$. We define {\it the marked Dynkin diagram of the variety of minimal rational tangents $\mathcal C_p\subset \mathbb P(T_pX)$} as {the marked Dynkin diagram of $Z$}.
\end{definition}

\begin{theorem}[Theorem \ref{theorem 2}]\label{thm: main 2}
      Let $X$ be the wonderful compactification of an irreducible adjoint symmetric homogeneous space $G/H$. Assume that the restricted root system is not of type A. Then, the marked Dynkin diagram of the variety of minimal rational tangents $\mathcal C_p \subset \mathbb P(T_pX)$ at a general point $p$ can be recovered from the marked Kac diagram of the open orbit $G/H$ by marking the adjacent node to a white node as Section \ref{Diagram}.
\end{theorem}

If $X$ is the wonderful compactification of the group type but not of type $A$, then the variety of minimal rational tangents at a general point $p \in X$ is isomorphic to $Z$ as Theorem \ref{Main} (1), see also Theorem 1.1 of \cite{BF15}. Let $P^{\alpha}$ be the maximal parabolic subgroup of $H$ corresponding to the marked roots $\alpha$ in the marked Dynkin diagram of $Z$. Then $Z=H/P^{\alpha}$ by Lemma \ref{isotoropy group type}.

If $X$ is the wonderful compactification of the simple type whose restricted root system is not of type $A$, then the variety of minimal rational tangents at a general point $p \in X$ is isomorphic to $Z$ by Theorem \ref{Main} (2), see also Proposition 3.14 of \cite{BKP}. If $G/H$ is not of type BI, CII, DI, then $H^0$ is simple and let $P^{\alpha}$ be the maximal parabolic subgroup of $H^0$ corresponding to the marked simple root $\alpha$ in the marked Dynkin diagram of $Z$. Then $Z=H^0/P^{\alpha}$ by Lemma \ref{isotoropy simple type}. If $G/H$ is not of type BI, CII, DI, then $H^0=(H_1^0\times H_2^0)/\Gamma$ decomposes into two simple factors and let $P^{\alpha_1}$ and $P^{\alpha_2}$ be the maximal parabolic subgroups of $H_1^0$ and $H_2^0$ corresponding to the two marked simple root $\alpha_1$ and $\alpha_2$ in the marked Dynkin diagram of $Z$. Then $Z=H^0/P^{\alpha_1}P^{\alpha_2}$ by Lemma \ref{isotoropy simple type}.

If $X$ is the wonderful compactification of the Hermitian exceptional type, then there are two isomorphic varieties of minimal rational tangents and two marked Dynkin diagram of $Z$. Assume that $G/H$ is not of type AIII. Let $P^{\alpha}$ be the maximal parabolic subgroup of the Levi subgroup $L$ corresponding to one marked simple root $\alpha$ in the one marked Dynkin diagram of $Z$. Then a variety of minimal rational tangents at a general point $p \in X$ is isomorphic to $L/P^{\alpha}$ by Theorem \ref{Main} (2) (see also Corollary 3.12 of \cite {BKP}) and by Lemma \ref{isotoropy hermitian type}. If $G/H$ is of type AIII, then the Levi subgroup $L$ decompose into two simple factors $L_1$ and $L_2$ and let $P^{\alpha_1}$ and $P^{\alpha_2}$ be the maximal parabolic subgroups of $L_1$ and $L_2$ corresponding to the two marked simple root $\alpha_1$ and $\alpha_2$ in one marked Dynkin diagram of $Z$. Then a variety of minimal rational tangents at a general point $p \in X$ is isomorphic to $Z=L/P^{\alpha_1}P^{\alpha_2}$ by Theorem \ref{Main} (2) (see also Corollary 3.12 of \cite {BKP}) and by Lemma \ref{isotoropy hermitian type}.

Let $X$ be the wonderful compactification of the Hermitian non-exceptional type not of type $A$. If $G/H$ is not of type AIII and let $P^{\alpha_1}$ and $P^{\alpha_2}$ be the two maximal parabolic subgroups of $L$ corresponding to two marked simple roots $\alpha_1$ and $\alpha_2$ respectively, in the marked Dynkin diagram of $Z$. Then the variety of minimal rational tangents at a general point $p \in X$ is isomorphic to the disjoint union of $L/P^{\alpha_1}$ and $L/P^{\alpha_2}$ by Theorem \ref{Main} (2) (see also Corollary 3.12 of \cite {BKP}) and by Lemma \ref{isotoropy hermitian type}. If $G/H$ is of type AIII, the marked Dynkin diagram of $Z$ consist of two Dynkin diagram of the Levi subgroup $L$ with four marked simple roots $\alpha_1$, $\alpha_2$, $\alpha_3$, $\alpha_4$. Since the Levi subgroup $L$ decompose into two simple factors $L_1$ and $L_2$, let $\alpha_1$ and $\alpha_2$ be the marked simple root of one Dynkin diagram of $L$ and let $P^{\alpha_1}$ and $P^{\alpha_2}$ be the maximal parabolic subgroups of $L_1$ and $L_2$ corresponding to the $\alpha_1$ and $\alpha_2$. Similarly, let $P^{\alpha_3}$ and $P^{\alpha_4}$ be the maximal parabolic subgroups of $L_1$ and $L_2$ corresponding to the $\alpha_3$ and $\alpha_4$. Then a variety of minimal rational tangents at a general point $p \in X$ is isomorphic to the disjoint union of $L/P^{\alpha_1}P^{\alpha_2}$ and $L/P^{\alpha_3}P^{\alpha_4}$ by Theorem \ref{Main} (2) (see also Corollary 3.12 of \cite {BKP}) and by Lemma \ref{isotoropy hermitian type}.

Hence, the marked Dynkin diagram of the variety of minimal rational tangents $\mathcal C_p\subset \mathbb P(T_pX)$ is well-defined except when the restricted root system is of type A.

\subsection{Marked Dynkin diagrams of VMRTs whose restricted root system is of type A}\label{A type}

We assume the restricted root system of $G/H$ is of type $A_r$, $r \geq 1$. Let $V_{\lambda}$ be the irreducible $G$-representation space with highest weight $\lambda$ such that $G/H$ is embedded as an open orbit in $\mathbb P(V_{\lambda})$ with a unique closed orbit. Due to \cite[Secion 4.5]{BKP}, we have possible pairs $(G, P_{\lambda})$ of a simple group $G$ up to finite covering and a parabolic subgroup $P_{\lambda}$ such that there exist an embedding $G/P_{\lambda}\subset\mathbb P(V_{\lambda})$ as the unique closed orbit. The list follows:
    \begin{center}
      \begin{tabular}{|c|c|c|c|c|}
        \hline
        $G/H$ & $A_r$ & $\lambda$ & $V_{\lambda}$ & $\frak p$ \\
        \hline
        $\PGL_{r+1}\times\PGL_{r+1}/\PGL_{r+1}$ & $r$ & $(\varpi_1,0) + (0,\varpi_r)$ 
        & ${\rm End}(\mathbb C^{r+1})$ & $\frak{sl}_{r+1}$ \\
        $\SL_{r+1}/{\rm SO}_{r+1}$ & $r$ & $2\varpi_1$ 
        & $S^2(\mathbb C^{r+1})$ & $S^2(\mathbb C^{r+1})_0$ \\
        $\SL_{2r+2}/\Sp_{2r+2}$ & $r$ & $\varpi_2$ 
        & $\Lambda^2(\mathbb C^{2r+2})$ & $\Lambda^2(\mathbb C^{2r+2})_0$ \\
        ${\rm SO}_{n}/{\rm S}(\OO_1 \times \OO_{n-1})$ & $1$ & $\varpi_1$ 
        & $\mathbb C^n$ & $\mathbb C^{n-1}$ \\
        $E_6/F_4$ & $2$ & $\varpi_1$ & $\mathbb C^{27}$ & $\mathbb C^{26}$ \\
         \hline
      \end{tabular}
    \end{center}
    
We recall the idea of proof of \cite[Theorem 4.41]{BKP}. There exist birational $G$-equivalent morphism $f:X\rightarrow \mathbb P(V_{\lambda})$ which is composition of the blow ups of strict transforms of orbit closures if $r\geq 2$. Because the set of general lines in $\mathbb P(V_{\lambda})$ is the hyperplane $\mathbb P(\frak p)$ and hence, the intersection of general lines with the closed orbit $G/P_{\lambda}$ in $\mathbb P(V_{\lambda})$ is the orbit $G/P_{\lambda}$ in $\mathbb P(\frak p)$.  By \cite[Proposition 9.7]{FH12}, the variety of minimal rational tangents $\mathcal C_p\subset \mathbb P(T_pX)$ is the intersection of general lines with the closed orbit $G/P_{\lambda}$, hence, the embedding $G/P_{\lambda}\subset\mathbb P(\frak p)$ is isomorphic to the variety of minimal rational tangents $\mathcal C_p\subset \mathbb P(T_pX)$.

\begin{lemma}\cite[Theorem 4.41]{BKP}\label{lem: A_r} Assume the restricted root system of $G/H$ is of type $A_r$, $r \geq 2$ and $X$ is the wonderful compactification of $G/H$. Then $G/H$ is either $\SL_{r+1}/\SO_{r+1}$, $(\PGL_{r+1}\times\PGL_{r+1})/\PGL_{r+1}$, $\SL_{2(r+1)}/\Sp_{2(r+1)}$, or $E_6/F_4$ and the variety of minimal rational tangents $\mathcal C_p\subset \mathbb P(T_pX)$ is either $\nu_2(\mathbb P^{r})$, $\mathbb P^{r} \times (\mathbb P^{r})^{\vee}$, $\Gr(2,2r+2)$, or $E_6/P_1$.\end{lemma}

We remark that the $H$-orbit $Z \subset \mathbb P \frak p$ is a general hyperplane section of $\mathcal C_p \subset \mathbb P(T_pX)$. On the other hand, if the restricted root system of $G/H$ is of type $A_1$, then $\mathcal C_p$ is isomorphic to the hyperplane $\mathbb P(\frak p)\subset\mathbb P(V_{\lambda})$, hence, $\mathcal C_p=\mathbb P(T_pX)=\mathbb P(\frak p)$.

\begin{lemma}\cite[Theorem 4.41]{BKP} Assume the restricted root system of $G/H$ is of type $A_1$ and $X$ is the wonderful compactification of $G/H$. Then $G/H$ is either $(\PGL_{2} \times \PGL_{2})/\PGL_{2}$, $\PGL_{2}/\PO_{2}$ or $\SO_{n+1}/\OO_{n}$ and the variety of minimal rational tangents $\mathcal C_p\subset \mathbb P(T_pX)$ is either $\mathbb P^2$, $\mathbb P^1$, or $\mathbb P^{n-1}$.
\end{lemma}

We have the marked Dynkin diagram of $(G, P_{\lambda})$, the Dynkin diagram with marked nodes corresponding to the parabolic subgroup $P_{\lambda}$.
\begin{definition}\label{vmrt L of A}
Let $X$ be the wonderful compactification of an irreducible adjoint symmetric homogeneous space $G/H$ whose restricted root system is of type $A$. %$A_r$, $r \geq 2$.
Let $\mathcal C_p \subset \mathbb P(T_pX)$ be a variety of minimal rational tangents at the base point $p \in X$. 
We define {\it the marked diagram of the variety of minimal rational tangents $\mathcal C_p\subset \mathbb P(T_pX)$} as (1) the marked Dynkin diagram of $(G, P_{\lambda})$ when the restricted root system is of type $A_r$ for $r \geq 2$. (2) the marked Dynkin diagram of the hyperplane $\mathbb P(\frak p)\subset\mathbb P(V_{\lambda})$ when the restricted root system is of type $A_1$.
\end{definition}
\begin{remark}
    To indicate the Veronese embedding $\nu_2(\mathbb P^{r}) \subset \mathbb P(\frak p)$, put $2$ on the crossed node of the marked diagram of the variety of minimal rational tangents $\mathcal C_p\subset \mathbb P(T_pX)$ when $G/H$ is $\SL_{r+1}/\SO_{r+1}$.
\end{remark}

Assume that $G/H$ is of the group type $\PGL_{r+1}\times\PGL_{r+1}/\PGL_{r+1}$. Then, $\sigma$ joins two Dyinkin diagrams of $H$ so that each node of the first diagram is joined with the respective node of the second diagram. This gives us a unique \emph{folding} from $H\times H$ to $H$. Hence, if $G/H$ is not of type $A_1$, the marked diagram of the variety of minimal rational tangents can be folded into the marked diagram of $Z \subset \mathbb P \frak p$. For example, 

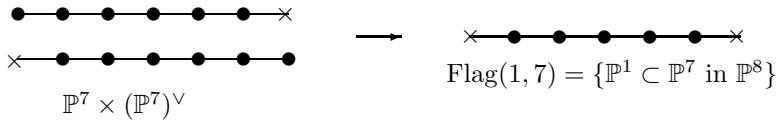
\setlength{\unitlength}{0.6cm}
\begin{center}
 \begin{picture}(0,3)(10,0)
\put(0,1.5){\line(1,0){1}}
\put(1,1.5){\line(1,0){.85}}
\put(2.15,1.5){\line(1,0){.9}}
\put(3,1.5){\line(1,0){.9}}
\put(4,1.5){\line(1,0){1}}
\put(5,1.5){\line(1,0){1}}

\put(0,1.5){\circle*{.3}}
\put(1,1.5){\circle*{.3}}
\put(2,1.5){\circle*{.3}}
\put(3,1.5){\circle*{.3}}
\put(4,1.5){\circle*{.3}}
\put(5,1.5){\circle*{.3}}
%\put(6,1.5){\circle*{.3}}

\put(0,0.5){\line(1,0){1}}
\put(1,0.5){\line(1,0){.85}}
\put(2.15,0.5){\line(1,0){.9}}
\put(3,0.5){\line(1,0){.9}}
\put(4,0.5){\line(1,0){1}}
\put(5,0.5){\line(1,0){1}}

%\put(0,0.5){\circle*{.3}}
\put(1,0.5){\circle*{.3}}
\put(2,0.5){\circle*{.3}}
\put(3,0.5){\circle*{.3}}
\put(4,0.5){\circle*{.3}}
\put(5,0.5){\circle*{.3}}
\put(6,0.5){\circle*{.3}}

\put(-0.3, 0.3){$\times$}
\put(5.7, 1.35){$\times$}

\put(7.5,1){\vector(1,0){1}}

\put(10,1){\line(1,0){1}}
\put(11,1){\line(1,0){1}}
\put(12,1){\line(1,0){1}}
\put(13,1){\line(1,0){1}}
\put(14,1){\line(1,0){1}}
\put(15,1){\line(1,0){1}}

\put(9.8, 0.85){$\times$}
\put(15.7, 0.85){$\times$}

\put(11,1){\circle*{.3}}
\put(12,1){\circle*{.3}}
\put(13,1){\circle*{.3}}
\put(14,1){\circle*{.3}}
\put(15,1){\circle*{.3}}

\put(1,-0.7){$\mathbb P^7 \times (\mathbb P^7)^{\vee} $}
\put(9.5,0){$\Flag(1,7)=\{ \mathbb P^1 \subset \mathbb P^7 \text{ in } \mathbb P^8 \}$ }
\end{picture}    
\end{center}
\vspace{1em}

Assume that $G/H$ is of a simple type. Then $G$ is either $\SL_{n}$, $\SO_n$, or $E_6$. The automorphism of the Dynkin diagram of $G$ is of order 2 which is the non-trivial element in the group of outer automorphisms $Aut(\frak g)/Int(\frak g)=\mathbb Z_2$. So we have $\mathbb Z_2$-gradation on $\frak g$
$$\frak g = \frak g_0 \oplus \frak g_1.$$
Hence, if $\frak g_0=\frak h$ and $\frak g_1 =\frak p$, then there exist a \emph{folding} from  Dynkin diagram of $G$ to  Dynkin diagram of $H$. Moreover, the pair $(G,H)$ is one of the followings: $(A_{2l}, B_l)$ for $l \geq 2$, $(A_{2l-1}, C_l)$ for $l \geq 3$, $(D_{l+1}, B_l)$ for $l \geq 2$, $(A_{2}, A_1)$, and $(E_6,F_4)$. We refer \cite[Section 8, Proposition 8.2]{Kac} for more details.

When $G/H$ is either $\SL_{2r+1}/\OO_{2r+1}$, $\SL_{2n}/\Sp_{2n}$ or $E_6/F_4$, there exist folding from $G$ to $H$ and the marked diagram of the variety of minimal rational tangents can be folded into the marked diagram of $Z \subset \mathbb P \frak p$, where the diagrams represent corresponding rational homogeneous spaces.
 Our examples are as follows;

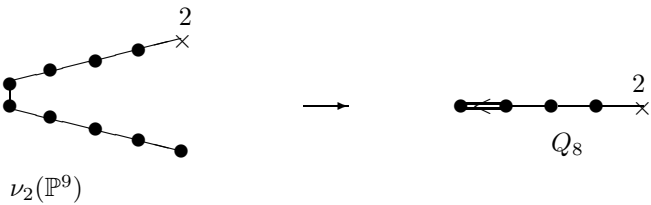
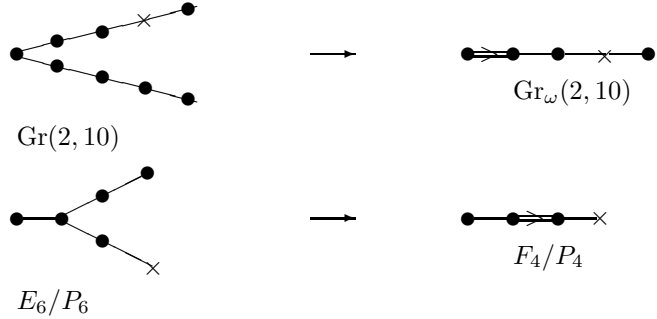
\setlength{\unitlength}{0.6cm}
\begin{center}
 \begin{picture}(0,4)(10,0)

\put(1,1.55){\line(4,1){3.8}}
\put(1,0.95){\line(0,1){0.5}}
\put(1,0.95){\line(4,-1){3.8}}

\put(1,1.5){\circle*{.3}}
\put(1,1){\circle*{.3}}
\put(1.9,1.8){\circle*{.3}}
\put(1.9,0.75){\circle*{.3}}
\put(2.9,2){\circle*{.3}}
\put(2.9,0.5){\circle*{.3}}
\put(4.6,2.3){$\times$}
\put(3.85,0.25){\circle*{.3}}
\put(3.85,2.25){\circle*{.3}}
\put(4.75,2.8){$2$}
\put(4.8,0){\circle*{.3}}

\put(7.5,1){\vector(1,0){1}}

\put(11,1.05){\line(1,0){1}}
\put(11.25,0.85){$<$}
\put(11,0.95){\line(1,0){1}}
\put(12,1){\line(1,0){.85}}
\put(13.15,1){\line(1,0){.9}}
\put(14.15,1){\line(1,0){.9}}

\put(14.8,1.3){$2$}
\put(14.8,0.8){$\times$}
\put(11,1){\circle*{.3}}
\put(12,1){\circle*{.3}}
\put(13,1){\circle*{.3}}
\put(14,1){\circle*{.3}}

\put(1,-1){$\nu_2(\mathbb P^9)$}
\put(13,0){$Q_{8}$}
\end{picture}
\vspace{1em}

\begin{picture}(0,3)(10,0)

\put(1,1.05){\line(4,1){4}}
\put(1,0.95){\line(4,-1){4}}

\put(1,1){\circle*{.3}}
\put(1.9,1.3){\circle*{.3}}
\put(1.9,0.75){\circle*{.3}}
\put(2.9,1.5){\circle*{.3}}
\put(2.9,0.5){\circle*{.3}}
\put(3.6,1.6){$\times$}
\put(3.85,0.25){\circle*{.3}}
\put(4.8,2){\circle*{.3}}
\put(4.8,0){\circle*{.3}}

\put(7.5,1){\vector(1,0){1}}

\put(11,1.05){\line(1,0){1}}
\put(11.25,0.85){$>$}
\put(11,0.95){\line(1,0){1}}
\put(12,1){\line(1,0){.85}}
\put(13.15,1){\line(1,0){.9}}
\put(14.15,1){\line(1,0){.9}}

\put(13.8,0.8){$\times$}
\put(11,1){\circle*{.3}}
\put(12,1){\circle*{.3}}
\put(13,1){\circle*{.3}}
\put(15,1){\circle*{.3}}

\put(1,-1){$\Gr(2,10)$}
\put(12,0){$\Gr_{\omega}(2,10)$}
\end{picture}
\vspace{1em}

\begin{picture}(0,3)(10,0)

\put(1.15,1){\line(1,0){0.9}}

\put(2,1.05){\line(2,1){2}}
\put(2,0.95){\line(2,-1){2}}

\put(1,1){\circle*{.3}}
\put(2,1){\circle*{.3}}
\put(2.9,1.5){\circle*{.3}}
\put(2.9,0.5){\circle*{.3}}
\put(3.9,2){\circle*{.3}}
\put(3.8,-.25){$\times$}

\put(7.5,1){\vector(1,0){1}}

\put(11.15,1){\line(1,0){0.9}}
\put(12,1.05){\line(1,0){1}}
\put(12.25,0.85){$>$}
\put(12,0.95){\line(1,0){1}}
\put(13,1){\line(1,0){.85}}

\put(11,1){\circle*{.3}}
\put(12,1){\circle*{.3}}
\put(13,1){\circle*{.3}}
\put(13.7,0.85){$\times$}

\put(1,-1){$E_6/P_6$}
\put(12,0){$F_4/P_4$}
\end{picture}
\vspace{2em}
\end{center}
Hence, we have the following.

\begin{theorem}[Theorem \ref{theorem 3}]\label{thm: main 3}
     Let $X$ be the wonderful compactification of an irreducible adjoint symmetric homogeneous space $G/H$. Assume that the restricted root system of $G/H$ is $A_r$ for $r \geq 2$. Then $G/H$ is either $\PGL_{r+1}\times\PGL_{r+1}/\PGL_{r+1}$, $\SL_{2r+1}/\OO_{2r+1}$, $\SL_{2n}/\Sp_{2n}$ or $E_6/F_4$ and the marked Dynkin diagram of the variety of minimal rational tangents $\mathcal C_p \subset \mathbb P(T_pX)$ at a general point $p$ can be recovered from the marked Kac diagram of the open orbit $G/H$  by marking the adjacent node to a white node as Section \ref{Diagram} and the folding process.
\end{theorem}

\begin{proof}
     If the Kac diagram of the open orbit $G/H$ is given, there is a marked diagram of $Z \subset \mathbb P \frak p$ as Section \ref{Diagram}. The marked diagram of $H$-orbit $Z \subset \mathbb P \frak p$ has a unique 2 to 1 covering diagram up to marking, which corresponds to the marked Dynkin diagram of the variety of minimal rational tangents.   
\end{proof}

\section{Appendix}

We give tables listing marked diagrams of VMRTs of adjoint irreducible symmetric spaces for each type. We follow \cite{bourbaki} the notion of types. We give marked Dynkin diagrams of the VMRT when the restricted root systems are of type $A$ in Table 4. In Table 5, we listing marked diagrams of VMRTs for group type when the restricted root system is not of type $A$. In Table 6, we are listing marked diagrams of VMRTs for simple types except of type \textrm{AI}, \textrm{AII}, \textrm{BII}, \textrm{DII}, and \textrm{EIV} in which the restricted root systems are of type $A$. We also give a table listing marked diagrams of VMRTs for Hermitian types except of type \textrm{AI} in Table 7.

\begin{table}
\centering
\begin{tabular}{|c|c| c | c |  c |}
\hline
type &$G$&$H$& VMRT & marked Dynkin diagram\\ \hline
Group type $A_1$ &$H\times H$ &$\PGL_{2}$ & $\mathbb P^{2}$ & 

\begin{picture}(2,1.5)(2,-0.3)
\put(2.25,-0.2){$\times$}
\put(2.5,0){\line(1,0){1}}
\put(3.5,0){\circle*{.3}}
\end{picture} \\ 

Group type $A_n$ &$H\times H$ &$\PGL_{n+1}$ & $\mathbb P^{n} \times (\mathbb P^{n})^{\vee} $ &

\begin{picture}(8,1.5)(-1,0.5)
\put(0,1.5){\line(1,0){1}}
\put(1,1.5){\line(1,0){.85}}
\put(2.1,1.5){\line(1,0){.4}}
\multiput(2.65,1.5)(0.2,0){5}{\line(1,0){.1}}
\put(3.55,1.5){\line(1,0){.4}}
\put(4,1.5){\line(1,0){1}}
\put(5,1.5){\line(1,0){1}}

\put(0,1.5){\circle*{.3}}
\put(1,1.5){\circle*{.3}}
\put(2,1.5){\circle*{.3}}
\put(4,1.5){\circle*{.3}}
\put(5,1.5){\circle*{.3}}

\put(0,0.5){\line(1,0){1}}
\put(1,0.5){\line(1,0){.85}}
\put(2.1,0.5){\line(1,0){.4}}
\multiput(2.65,0.5)(0.2,0){5}{\line(1,0){.1}}
\put(3.55,0.5){\line(1,0){.4}}
\put(4,0.5){\line(1,0){1}}
\put(5,0.5){\line(1,0){1}}

\put(1,0.5){\circle*{.3}}
\put(2,0.5){\circle*{.3}}
\put(4,0.5){\circle*{.3}}
\put(5,0.5){\circle*{.3}}
\put(6,0.5){\circle*{.3}}

\put(-0.3, 0.3){$\times$}
\put(5.7, 1.35){$\times$}
\end{picture}
\\

Herm. non-excep.AI &$\PGL_{2}$&$\PO_{2}$ & $\mathbb P^{1}$&
\begin{picture}(2,1.5)(2,-0.3)
\put(2.75,-0.2){$\times$}
\end{picture}
\\ 

AI &$\SL_{n}$&$\SO_{n}$ & $\nu_2(\mathbb P^{n-1})$&
\begin{picture}(8,1.5)(-1,-0.3)
\put(0,0){\line(1,0){1}}
\put(1,0){\line(1,0){.85}}
\put(2.1,0){\line(1,0){.4}}
\multiput(2.65,0)(0.2,0){5}{\line(1,0){.1}}
\put(3.55,0){\line(1,0){.4}}
\put(4,0){\line(1,0){1}}
\put(5,0){\line(1,0){1}}

\put(1,0){\circle*{.3}}
\put(2,0){\circle*{.3}}
\put(4,0){\circle*{.3}}
\put(5,0){\circle*{.3}}
\put(6,0){\circle*{.3}}

\put(-.15, -.2){$\times$}
\put(-.10, 0.2){$2$}

\end{picture} \\

AII &$\SL_{2n}$&$\Sp_{2n}$ &$\Gr(2,2n)$&
\begin{picture}(3,1.5)(1.5,0.2)
\put(0,0.5){\line(1,0){1}}
\put(1,0.5){\line(1,0){.85}}
\put(2.1,0.5){\line(1,0){.4}}
\multiput(2.65,0.5)(0.2,0){5}{\line(1,0){.1}}
\put(3.55,0.5){\line(1,0){.4}}
\put(4,0.5){\line(1,0){1}}
\put(5,0.5){\line(1,0){1}}

\put(0,0.5){\circle*{.3}}
\put(2,0.5){\circle*{.3}}
\put(4,0.5){\circle*{.3}}
\put(5,0.5){\circle*{.3}}
\put(6,0.5){\circle*{.3}}
\put(0.7, 0.3){$\times$}
\end{picture}
\\

BII & $\SO_{2n+1}$&$\OO_{2n}$ & $\mathbb P^{2n-1}$&

\begin{picture}(2,1.5)(2,-0.3)

\put(0,0){\line(1,0){1}}
\put(1,0){\line(1,0){.85}}
\put(2.1,0){\line(1,0){.4}}
\multiput(2.65,0)(0.2,0){5}{\line(1,0){.1}}
\put(3.55,0){\line(1,0){.4}}
\put(4,0){\line(1,0){1}}
\put(5,0){\line(1,0){1}}
\put(1,0){\circle*{.3}}
\put(2,0){\circle*{.3}}
\put(4,0){\circle*{.3}}
\put(5,0){\circle*{.3}}
\put(6,0){\circle*{.3}}

\put(-.25, -.2){$\times$}
\end{picture}
\\

DII &$\SO_{2n}$&$\OO_{2n-1}$ &  $\mathbb P^{2n-2}$ &

\begin{picture}(2,1.5)(2,-0.2)
\put(0,0){\line(1,0){1}}
\put(1,0){\line(1,0){.85}}
\put(2.1,0){\line(1,0){.4}}
\multiput(2.65,0)(0.2,0){5}{\line(1,0){.1}}
\put(3.55,0){\line(1,0){.4}}
\put(4,0){\line(1,0){1}}
\put(5,0){\line(1,0){1}}

\put(1,0){\circle*{.3}}
\put(2,0){\circle*{.3}}
\put(4,0){\circle*{.3}}
\put(5,0){\circle*{.3}}
\put(6,0){\circle*{.3}}
\put(-.25, -.2){$\times$}
\end{picture}
\\

EIV &$E_6$&$F_4$ &$E_6/P_1$ &

\begin{picture}(2,1.5)(2,-0.5)
\put(1,0){\line(1,0){1}}
\put(2,0){\line(1,0){1}}
\put(3,0){\line(1,0){1}}
\put(4,0){\line(1,0){1}}
\put(3,0){\line(0,-1){0.8}}
\put(2,0){\circle*{.3}}
\put(3,0){\circle*{.3}}
\put(4,0){\circle*{.3}}
\put(5,0){\circle*{.3}}
\put(3,-1){\circle*{.3}}
\put(0.75,-.15){$\times$}
\end{picture}
\\ &&&& \\\hline
\end{tabular}\\
\caption{The restricted root system is of type A}
\end{table}

\begin{table}
\centering
    \begin{tabular}{|c|c|c|c|} \hline
type & $H$ & VMRT & marked Dynkin diagram\\ \hline
&&&\\
$A_n$&$\PGL_{n+1}$ &  & \\

$B_n$ & $\PO_{2n+1}$ &$\OG(2,2n+1)$ & 

\begin{picture}(2,1.5)(2,-0.3)
\put(0,0){\line(1,0){1}}
\put(1,0){\line(1,0){.85}}
\put(2.1,0){\line(1,0){.4}}
\multiput(2.65,0)(0.2,0){5}{\line(1,0){.1}}
\put(3.55,0){\line(1,0){.4}}
\put(4,0){\line(1,0){1}}
\put(5,0.05){\line(1,0){1}}
\put(5,-0.05){\line(1,0){1}}

\put(0,0){\circle*{.3}}
\put(2,0){\circle*{.3}}
\put(4,0){\circle*{.3}}
\put(5,0){\circle*{.3}}
\put(6,0){\circle*{.3}}

\put(0.75, -.2){$\times$}
\put(5.2, -.15){$>$}
\end{picture}
\\

$C_n$& $\PSp_{2n}$ & $v_2(\mathbb P^{2n-1})$ &

\begin{picture}(8,1.5)(-1,-0.3)
\put(0,0){\line(1,0){1}}
\put(1,0){\line(1,0){.85}}
\put(2.1,0){\line(1,0){.4}}
\multiput(2.65,0)(0.2,0){5}{\line(1,0){.1}}
\put(3.55,0){\line(1,0){.4}}
\put(4,0){\line(1,0){1}}
\put(5,0.05){\line(1,0){1}}
\put(5,-0.05){\line(1,0){1}}

\put(1,0){\circle*{.3}}
\put(2,0){\circle*{.3}}
\put(4,0){\circle*{.3}}
\put(5,0){\circle*{.3}}
\put(6,0){\circle*{.3}}

\put(0, 0.2){$2$}
\put(-.15, -.2){$\times$}
\put(5.2, -.15){$<$}
\end{picture}
\\

$D_n$& $\PO_{2n}$ & $\OG(2,2n)$ &
\begin{picture}(2,1.5)(2,-0.2)
\put(0,0){\line(1,0){1}}
\put(1,0){\line(1,0){.85}}
\put(2.1,0){\line(1,0){.4}}
\multiput(2.65,0)(0.2,0){5}{\line(1,0){.1}}
\put(3.55,0){\line(1,0){.4}}
\put(4,0){\line(1,0){1}}
\put(5,0){\line(1,1){0.7}}
\put(5,0){\line(1,-1){0.7}}

\put(0,0){\circle*{.3}}
\put(2,0){\circle*{.3}}
\put(4,0){\circle*{.3}}
\put(5,0){\circle*{.3}}
\put(5.7,0.65){\circle*{.3}}
\put(5.7,-0.65){\circle*{.3}}
\put(.75, -.2){$\times$}
\end{picture}
\\
$E_6$ & $E^{ad}_6$ & $E_6/P_2$ & 

\begin{picture}(2,1.5)(3,-0.5)
%\put(0,0){\line(1,0){1}}
\put(1,0){\line(1,0){1}}
\put(2,0){\line(1,0){1}}
\put(3,0){\line(1,0){1}}
\put(4,0){\line(1,0){1}}

\put(3,0){\line(0,-1){0.8}}

\put(1,0){\circle*{.3}}
\put(2,0){\circle*{.3}}
\put(3,0){\circle*{.3}}
\put(4,0){\circle*{.3}}
\put(5,0){\circle*{.3}}
%\put(6,0){\circle*{.3}}
\put(2.75,-1){$\times$}
\end{picture}
 \\
$E_7$ & $E^{ad}_7$ & $E_7/P_1$ & 

\begin{picture}(2,1.5)(3,-0.5)
\put(1,0){\line(1,0){1}}
\put(2,0){\line(1,0){1}}
\put(3,0){\line(1,0){1}}
\put(4,0){\line(1,0){1}}
\put(5,0){\line(1,0){1}}

\put(3,0){\line(0,-1){0.8}}

\put(2,0){\circle*{.3}}
\put(3,0){\circle*{.3}}
\put(4,0){\circle*{.3}}
\put(5,0){\circle*{.3}}
\put(6,0){\circle*{.3}}
\put(3,-1){\circle*{.3}}
\put(0.75,-.15){$\times$}
\end{picture}

\\
$E_8$ & $E^{ad}_8$ & $E_8/P_8$ & 

\begin{picture}(2,1.5)(2,-0.5)
\put(0,0){\line(1,0){1}}
\put(1,0){\line(1,0){1}}
\put(2,0){\line(1,0){1}}
\put(3,0){\line(1,0){1}}
\put(4,0){\line(1,0){1}}
\put(5,0){\line(1,0){1}}

\put(2,0){\line(0,-1){0.8}}

\put(0,0){\circle*{.3}}
\put(1,0){\circle*{.3}}
\put(2,0){\circle*{.3}}
\put(3,0){\circle*{.3}}
\put(4,0){\circle*{.3}}
\put(5,0){\circle*{.3}}
%\put(6,0){\circle*{.3}}

\put(2,-1){\circle*{.3}}
\put(5.75,-.15){$\times$}
\end{picture}
 \\

$F_4$ & $F^{ad}_4$ &$F_4/P_1$ & 

\begin{picture}(2,1.5)(2,0)
\put(1,0){\line(1,0){1}}
\put(3,0){\line(1,0){1}}
\put(2,0.05){\line(1,0){1}}
\put(2,-0.05){\line(1,0){1}}

\put(2,0){\circle*{.3}}
\put(3,0){\circle*{.3}}
\put(4,0){\circle*{.3}}

\put(0.75, -.2){$\times$}
\put(2.2, -.15){$>$}
\end{picture}
\\
$G_2$ & $G^{ad}_2$ &$G_2/P_2$ & 

 \begin{picture}(4,1.5)(1,0)
\put(2,0){\line(1,0){1}}
\put(2,0.07){\line(1,0){1}}
\put(2,-0.07){\line(1,0){1}}

\put(2,0){\circle*{.3}}
\put(2.75, -.2){$\times$}
\put(2.2, -.15){$<$}
\end{picture}   
\\
&&& \\
\hline
  \end{tabular}
\caption{Group type $(H\times H)/H$}
\end{table}

\begin{table}
\centering
\begin{tabular}{|c|c|c|c|c|c|} 
\hline
type& $G$ & $H$ & VMRT & marked Dynkin diagram \\ 
\hline
&&&&\\
AI &$\SL_{n}$&$\SO_{n}$ & & \\

AII &$\SL_{2n}$&$\Sp_{2n}$ & & \\

\textrm{BI}& $\SO_{2n+1}$ & $\rm S(\OO_m\times \OO_{2n+1-m})$& $Q_m\times Q_{2n+1-m}$ &

\begin{picture}(8,1.5)(-0.5,0)
\put(0.3,1.15){\line(1,-1){0.7}}
\put(0.3,-0.15){\line(1,1){0.7}}
\multiput(1,0.5)(0.2,0){5}{\line(1,0){.1}}
\put(2,0.5){\line(1,0){1}}
\put(4,0.5){\line(1,0){1}}
\multiput(5,0.5)(0.2,0){5}{\line(1,0){.1}}
\put(6,0.45){\line(1,0){1}}
\put(6,0.55){\line(1,0){1}}

\put(0.3,1.15){\circle*{.3}}
\put(0.3,-0.15){\circle*{.3}}
\put(1,0.5){\circle*{.3}}
\put(2,0.5){\circle*{.3}}
\put(5,0.5){\circle*{.3}}
\put(6,0.5){\circle*{.3}}
\put(7,0.5){\circle*{.3}}
\put(2.7, 0.35){$\times$}
\put(3.7, 0.35){$\times$}
\put(6.25,0.35){$>$}
\end{picture}\\

BII & $\SO_{2n+1}$&$\OO_{2n}$ & &\\

\textrm{CII}& $\Sp_{2n}$ & $\Sp_{2m}\times \Sp_{2(n-m)}$& $\mathbb P^{2m-1}\times \mathbb P^{2(n-m)-1}$ &

\begin{picture}(8,1.5)(-0.5,0)
\put(0,0.45){\line(1,0){1}}
\put(0,0.55){\line(1,0){1}}
\multiput(1,0.5)(0.2,0){5}{\line(1,0){.1}}
\put(2,0.5){\line(1,0){1}}
\put(4,0.5){\line(1,0){1}}
\multiput(5,0.5)(0.2,0){5}{\line(1,0){.1}}
\put(6,0.45){\line(1,0){1}}
\put(6,0.55){\line(1,0){1}}

\put(1,0.5){\circle*{.3}}
\put(2,0.5){\circle*{.3}}
\put(0,0.5){\circle*{.3}}
\put(5,0.5){\circle*{.3}}
\put(6,0.5){\circle*{.3}}
\put(7,0.5){\circle*{.3}}
\put(2.7, 0.35){$\times$}
\put(3.7, 0.35){$\times$}
\put(6.3,0.35){$<$}
\put(0.25,0.35){$>$}
\end{picture}
\\

\textrm{DI}& $\SO_{2n}$ &  $\rm S(\OO_{2m+1}\times \OO_{2(n-m)-1})$& $Q_{2m+1}\times Q_{2(n-m)-1}$ &

\begin{picture}(8,1.5)(-0.5,0)
\put(0,0.45){\line(1,0){1}}
\put(0,0.55){\line(1,0){1}}
\multiput(1,0.5)(0.2,0){5}{\line(1,0){.1}}
\put(2,0.5){\line(1,0){1}}
\put(4,0.5){\line(1,0){1}}
\multiput(5,0.5)(0.2,0){5}{\line(1,0){.1}}
\put(6,0.45){\line(1,0){1}}
\put(6,0.55){\line(1,0){1}}

\put(1,0.5){\circle*{.3}}
\put(2,0.5){\circle*{.3}}
\put(0,0.5){\circle*{.3}}
\put(5,0.5){\circle*{.3}}
\put(6,0.5){\circle*{.3}}
\put(7,0.5){\circle*{.3}}
\put(2.7, 0.35){$\times$}
\put(3.7, 0.35){$\times$}
\put(6.3,0.35){$>$}
\put(0.25,0.35){$<$}
\end{picture} \\

& &  $\rm S(\OO_{2m}\times \OO_{2(n-m)})$& $Q_{2m}\times Q_{2(n-m)}$ &

\begin{picture}(8,2)(-0.5,0)
\put(0.3,1.15){\line(1,-1){0.7}}
\put(0.3,-0.15){\line(1,1){0.7}}
\multiput(1,0.5)(0.2,0){5}{\line(1,0){.1}}
\put(2,0.5){\line(1,0){1}}
\put(4,0.5){\line(1,0){1}}
\multiput(5,0.5)(0.2,0){5}{\line(1,0){.1}}
\put(5,0.5){\line(1,0){1}}
\put(6,0.5){\line(1,1){0.7}}
\put(6,0.5){\line(1,-1){0.7}}

\put(0.3,1.15){\circle*{.3}}
\put(0.3,-0.15){\circle*{.3}}
\put(1,0.5){\circle*{.3}}
\put(2,0.5){\circle*{.3}}
\put(5,0.5){\circle*{.3}}
\put(6,0.5){\circle*{.3}}
\put(6.7,1.15){\circle*{.3}}
\put(6.7,-0.15){\circle*{.3}}

\put(2.7, 0.35){$\times$}
\put(3.7, 0.35){$\times$}
\end{picture}
\\

DII &$\SO_{2n}$&$\OO_{2n-1}$ & & \\

\textrm{EI}& $E_6$ & $C_4$& $\LG(4,8)$ &

\begin{picture}(8,1.5)(-1,0)
\put(1,0){\line(1,0){1}}
\put(2,0){\line(1,0){1}}
\put(3,0.05){\line(1,0){1}}
\put(3,-0.05){\line(1,0){1}}
\put(1,0){\circle*{.3}}
\put(2,0){\circle*{.3}}
\put(3,0){\circle*{.3}}
\put(3.75,-.15){$\times$}
\put(3.25,-.15){$<$}
\end{picture}
 \\
\textrm{EII}& $E_6$ & $A_5 \times A_1$& $\Gr(3,6)\times \mathbb P^1$ &

\begin{picture}(8,1.5)(0,0)
\put(1,0){\line(1,0){1}}
\put(2,0){\line(1,0){1}}
\put(3,0){\line(1,0){1}}
\put(4,0){\line(1,0){1}}
\put(2,0){\circle*{.3}}
\put(1,0){\circle*{.3}}
\put(4,0){\circle*{.3}}
\put(5,0){\circle*{.3}}
\put(2.75,-.16){$\times$}
\put(5.75,-.16){$\times$}
\end{picture}
\\

EIV &$E_6$&$F_4$ & & \\

\textrm{EV}& $E_7$ & $A_7$ & $\Gr(4,8)$ &

\begin{picture}(8,1.5)(-1,1)
\put(0,1){\line(1,0){1}}
\put(1,1){\line(1,0){.85}}
\put(2.15,1){\line(1,0){.9}}
\put(3,1){\line(1,0){.9}}
\put(4,1){\line(1,0){1}}
\put(5,1){\line(1,0){1}}

\put(0,1){\circle*{.3}}
\put(1,1){\circle*{.3}}
\put(2,1){\circle*{.3}}
%\put(3,1){\circle*{.3}}
\put(4,1){\circle*{.3}}
\put(5,1){\circle*{.3}}
\put(6,1){\circle*{.3}}
\put(2.75,0.84){$\times$}
\end{picture}
 \\
\textrm{EVI}& $E_7$ & $D_6 \times A_1$ & $\OG(6,12)\times \mathbb P^1$ &

\begin{picture}(2,1.5)(0,0)
\put(-2,0){\line(1,0){1}}
\put(-1,0){\line(1,0){1}}
\put(0,0){\line(1,0){1}}
\put(1,0){\line(1,1){0.7}}
\put(1,0){\line(1,-1){0.7}}
\put(-2,0){\circle*{.3}}
\put(-1,0){\circle*{.3}}
\put(0,0){\circle*{.3}}
\put(1,0){\circle*{.3}}
\put(1.7,-0.65){\circle*{.3}}

\put(3, -0.16){$\times$}
\put(1.4, 0.53){$\times$}
\end{picture}
\\
\textrm{EVIII}& $E_8$ & $D_8$& $\OG(8,16)$ & 

\begin{picture}(2,2)(1,-0.5)
\put(-1,0){\line(1,0){1}}
\put(0,0){\line(1,0){1}}
\put(1,0){\line(1,0){1}}
\put(2,0){\line(1,0){1}}
\put(3,0){\line(1,0){1}}
\put(4,0){\line(1,1){0.7}}
\put(4,0){\line(1,-1){0.7}}
\put(-1,0){\circle*{.3}}
\put(0,0){\circle*{.3}}
\put(1,0){\circle*{.3}}
\put(2,0){\circle*{.3}}
\put(3,0){\circle*{.3}}
\put(4,0){\circle*{.3}}
\put(4.7,-0.65){\circle*{.3}}
\put(4.4, 0.53){$\times$}
\end{picture}

\\
\textrm{EVIII}& $E_8$ & $E_7\times A_1$& $E_7/P_7\times \mathbb P^1$ &

\begin{picture}(7,2)(-0.5,0)
\put(0,0){\line(1,0){1}}
\put(1,0){\line(1,0){1}}
\put(2,0){\line(1,0){1}}
\put(3,0){\line(1,0){1}}
\put(4,0){\line(1,0){1}}
\put(2,0){\line(0,-1){0.8}}
\put(0,0){\circle*{.3}}
\put(1,0){\circle*{.3}}
\put(2,0){\circle*{.3}}
\put(3,0){\circle*{.3}}
\put(4,0){\circle*{.3}}
\put(2,-1){\circle*{.3}}
\put(4.75,-.15){$\times$}
\put(5.75,-.15){$\times$}
\end{picture}\\
\textrm{FI}& $F_4$ & $C_3\times A_1$& $\LG(3,6)\times\mathbb P^1$ &

\begin{picture}(0,2)(4,0)
\put(2,0){\line(1,0){1}}
\put(3,0.05){\line(1,0){1}}
\put(3,-0.05){\line(1,0){1}}
\put(2,0){\circle*{.3}}
\put(3,0){\circle*{.3}}
\put(3.75,-.15){$\times$}
\put(3.25,-.15){$<$}
\put(4.75,-.16){$\times$}
\end{picture}
 \\
\textrm{FII}& $F_4$ & $B_4$& $\OG(4,9)$ &

\begin{picture}(2,1.5)(2,0)
\put(1,0){\line(1,0){1}}
\put(2,0){\line(1,0){1}}
\put(3,0.05){\line(1,0){1}}
\put(3,-0.05){\line(1,0){1}}
\put(1,0){\circle*{.3}}
\put(2,0){\circle*{.3}}
\put(3,0){\circle*{.3}}
\put(3.75,-.15){$\times$}
\put(3.25,-.15){$>$}
\end{picture}\\
\textrm{G}& $G_2$ & $A_1\times A_1$ & $\mathbb P^1 \times \nu_3(\mathbb P^1)$ &

\begin{picture}(6,1.5)(0,0)
\put(1.7,0){$\times$}
\put(2.7,0){$\times$}
\put(2.8, 0.4){$3$}
\end{picture}
\\&&&&\\ \hline
  \end{tabular} \\
\caption {Simple type $G/H$}
\end{table}

\begin{table}
\centering
\begin{tabular}{|c|c|c|c|c|}
\hline
type& $G$ & $L$ & VMRT & marked Dynkin diagram\\ 
\hline
\textrm{AIII}&$\PGL_{n}$ & $\rm PG(L_{m}\times L_{n-m})$& $ \mathbb P^{m-1} \times (\mathbb P^{n-m-1})^{\vee}$ &
\begin{picture}(10,1.5)(-1,1.5)
\put(0,1.5){\line(1,0){1}}
\multiput(1,1.5)(0.2,0){5}{\line(1,0){.1}}
\put(2,1.5){\line(1,0){1}}
\put(4,1.5){\line(1,0){1}}
\multiput(5,1.5)(0.2,0){5}{\line(1,0){.1}}
\put(6,1.5){\line(1,0){1}}
\put(7,1.5){\line(1,0){1}}

\put(1,1.5){\circle*{.3}}
\put(2,1.5){\circle*{.3}}
\put(3,1.5){\circle*{.3}}
\put(4,1.5){\circle*{.3}}
\put(5,1.5){\circle*{.3}}
\put(6,1.5){\circle*{.3}}
\put(7,1.5){\circle*{.3}}
\put(-0.3, 1.35){$\times$}
\put(7.7, 1.35){$\times$}

\put(0,0.5){\line(1,0){1}}
\multiput(1,0.5)(0.2,0){5}{\line(1,0){.1}}
\put(2,0.5){\line(1,0){1}}
\put(4,0.5){\line(1,0){1}}
\multiput(5,0.5)(0.2,0){5}{\line(1,0){.1}}
\put(6,0.5){\line(1,0){1}}
\put(7,0.5){\line(1,0){1}}
\put(1,0.5){\circle*{.3}}
\put(2,0.5){\circle*{.3}}
\put(0,0.5){\circle*{.3}}
\put(7,0.5){\circle*{.3}}
\put(5,0.5){\circle*{.3}}
\put(6,0.5){\circle*{.3}}
\put(8,0.5){\circle*{.3}}
\put(2.7, 0.35){$\times$}
\put(3.7, 0.35){$\times$}
\end{picture}
\\
&&& and $(\mathbb P^{m-1})^{\vee}\times\mathbb P^{n-m-1} $&\\

\textrm{DIII}& $\PO_{4n+2}$& $\PGL_{2n+1}$ &$\Gr(2n-1,2n+1)$ & 

\begin{picture}(8,2)(0,1.5)
\put(0,0.5){\line(1,0){1}}
\put(1,0.5){\line(1,0){.85}}
\put(2.1,0.5){\line(1,0){.4}}
\multiput(2.65,0.5)(0.2,0){5}{\line(1,0){.1}}
\put(3.55,0.5){\line(1,0){.4}}
\put(4,0.5){\line(1,0){1}}
\put(5,0.5){\line(1,0){1}}

\put(0,0.5){\circle*{.3}}
\put(2,0.5){\circle*{.3}}
\put(4,0.5){\circle*{.3}}
\put(5,0.5){\circle*{.3}}
\put(6,0.5){\circle*{.3}}
\put(0.7, 0.3){$\times$}

\put(0,1.5){\line(1,0){1}}
\put(1,1.5){\line(1,0){.85}}
\put(2.1,1.5){\line(1,0){.4}}
\multiput(2.65,1.5)(0.2,0){5}{\line(1,0){.1}}
\put(3.55,1.5){\line(1,0){.4}}
\put(4,1.5){\line(1,0){1}}
\put(5,1.5){\line(1,0){1}}

\put(0,1.5){\circle*{.3}}
\put(2,1.5){\circle*{.3}}
\put(4,1.5){\circle*{.3}}
\put(1,1.5){\circle*{.3}}
\put(6,1.5){\circle*{.3}}
\put(4.7, 1.3){$\times$}
%\put(6.5, 0.8){$\updownarrow$$\sigma$}
\end{picture}
\\
&&& and $\Gr(2,2n+1)$&\\

\textrm{EIII}& $E^{ad}_6$ & $D_5 \times \mathbb C^*$ & $\OG(5,10)$ &
\begin{picture}(2,2)(3,0)
\put(0,0){\line(1,0){1}}
\put(1,0){\line(1,0){1}}
\put(2,0){\line(1,1){0.7}}
\put(2,0){\line(1,-1){0.7}}

\put(0,0){\circle*{.3}}
\put(1,0){\circle*{.3}}
\put(2,0){\circle*{.3}}
%\put(1.7,0.65){\circle*{.3}}
\put(2.7,-0.65){\circle*{.3}}

\put(4,0){\line(1,0){1}}
\put(5,0){\line(1,0){1}}
\put(6,0){\line(1,1){0.7}}
\put(6,0){\line(1,-1){0.7}}

\put(4,0){\circle*{.3}}
\put(5,0){\circle*{.3}}
\put(6,0){\circle*{.3}}
\put(6.7,0.65){\circle*{.3}}
%\put(5.7,-0.65){\circle*{.3}}

\put(6.4, -0.85){$\times$}
\put(2.4, 0.53){$\times$}
\end{picture}
\\
&&&and $\OG(5,10)$ & \\
&&&&  \\
\hline 
\hline

\textrm{AI}&$\PGL_{2}$ &$\PO_{2}$ & & \\

\textrm{AIII}&$\PGL_{2n}$ & $\rm PG(L_{n}\times L_n)$ & 
$\mathbb P^{n-1} \times (\mathbb P^{n-1})^{\vee}$&

\begin{picture}(6,1.5)(1,1.5)
\put(0,2){\line(1,0){1}}
\multiput(1,2)(0.2,0){5}{\line(1,0){.1}}
\put(2,2){\line(1,0){1}}
\put(4,2){\line(1,0){1}}
\multiput(5,2)(0.2,0){5}{\line(1,0){.1}}
\put(6,2){\line(1,0){1}}

\put(1,2){\circle*{.3}}
\put(2,2){\circle*{.3}}
\put(3,2){\circle*{.3}}
\put(4,2){\circle*{.3}}
\put(5,2){\circle*{.3}}
\put(6,2){\circle*{.3}}

\put(-0.3, 1.85){$\times$}
\put(6.7, 1.85){$\times$}

\put(0,1){\line(1,0){1}}
\multiput(1,1)(0.2,0){5}{\line(1,0){.1}}
\put(2,1){\line(1,0){1}}
\put(4,1){\line(1,0){1}}
\multiput(5,1)(0.2,0){5}{\line(1,0){.1}}
\put(6,1){\line(1,0){1}}

\put(1,1){\circle*{.3}}
\put(2,1){\circle*{.3}}
\put(0,1){\circle*{.3}}
\put(7,1){\circle*{.3}}
\put(5,1){\circle*{.3}}
\put(6,1){\circle*{.3}}

\put(3.7, 0.85){$\times$}
\put(2.7, 0.85){$\times$}
\put(7.5, 1.3){$\updownarrow$$\sigma$}
\end{picture}\\
&&& $\sqcup (\mathbb P^{n-1})^{\vee} \times \mathbb P^{n-1}$ & \\

\textrm{BI}& $\PO_{2n+1}$&$\rm P(\OO_2 \times \OO_{2n-1})$ &$Q_{2n-3}\sqcup Q_{2n-3}$ & 

\begin{picture}(8,2)(0,0.5)
\put(0,0){\line(1,0){1}}
\put(1,0){\line(1,0){.85}}
\put(2.1,0){\line(1,0){.4}}
\multiput(2.65,0)(0.2,0){5}{\line(1,0){.1}}
\put(3.55,0){\line(1,0){.4}}
\put(4,0){\line(1,0){1}}
\put(5,0.05){\line(1,0){1}}
\put(5,-0.05){\line(1,0){1}}

\put(1,0){\circle*{.3}}
\put(2,0){\circle*{.3}}
\put(4,0){\circle*{.3}}
\put(5,0){\circle*{.3}}
\put(6,0){\circle*{.3}}

\put(-.25, -.2){$\times$}
\put(5.2, -.15){$>$}

\put(0,1){\line(1,0){1}}
\put(1,1){\line(1,0){.85}}
\put(2.1,1){\line(1,0){.4}}
\multiput(2.65,1)(0.2,0){5}{\line(1,0){.1}}
\put(3.55,1){\line(1,0){.4}}
\put(4,1){\line(1,0){1}}
\put(5,1.05){\line(1,0){1}}
\put(5,0.95){\line(1,0){1}}

\put(1,1){\circle*{.3}}
\put(2,1){\circle*{.3}}
\put(4,1){\circle*{.3}}
\put(5,1){\circle*{.3}}
\put(6,1){\circle*{.3}}

\put(-.25, 0.8){$\times$}
\put(5.2, 0.85){$>$}
\put(6.5, 0.25){$\updownarrow$$\sigma$}
\end{picture} 
\\

\textrm{CI}& $\PSp_{2n}$&$\PGL_n$ & $v_2(\mathbb P^{n-1}) \sqcup v_2(\mathbb P^{n-1})^{\vee}$ &
\begin{picture}(8,3)(0,0.5)
\put(0,0){\line(1,0){1}}
\put(1,0){\line(1,0){.85}}
\put(2.1,0){\line(1,0){.4}}
\multiput(2.65,0)(0.2,0){5}{\line(1,0){.1}}
\put(3.55,0){\line(1,0){.4}}
\put(4,0){\line(1,0){1}}
\put(5,0){\line(1,0){1}}

\put(1,0){\circle*{.3}}
\put(2,0){\circle*{.3}}
\put(4,0){\circle*{.3}}
\put(5,0){\circle*{.3}}
\put(6,0){\circle*{.3}}

\put(-0.2, 0.2){$2$}
\put(-0.3, -0.2){$\times$}

\put(0,1){\line(1,0){1}}
\put(1,1){\line(1,0){.85}}
\put(2.1,1){\line(1,0){.4}}
\multiput(2.65,1)(0.2,0){5}{\line(1,0){.1}}
\put(3.55,1){\line(1,0){.4}}
\put(4,1){\line(1,0){1}}
\put(5,1){\line(1,0){1}}

\put(0,1){\circle*{.3}}
\put(1,1){\circle*{.3}}
\put(2,1){\circle*{.3}}
\put(4,1){\circle*{.3}}
\put(5,1){\circle*{.3}}

\put(58, 1.2){$2$}
\put(5.8, 1.2){$2$}
\put(5.7, 0.8){$\times$}
\put(6.5, 0.3){$\updownarrow$$\sigma$}
\end{picture}
\\ &&&& \\

\textrm{DI}& $\PO_{2n}$&$\rm P(\OO_2 \times \OO_{2(n-1)})$ &$Q_{2n-4} \sqcup Q_{2n-4}$ & 

\begin{picture}(8,3)(0,0.5)
\put(0,0){\line(1,0){1}}
\put(1,0){\line(1,0){.85}}
\put(2.1,0){\line(1,0){.4}}
\multiput(2.65,0)(0.2,0){5}{\line(1,0){.1}}
\put(3.55,0){\line(1,0){.4}}
\put(4,0){\line(1,0){1}}
\put(5,0){\line(1,1){0.6}}
\put(5,0){\line(1,-1){0.6}}

\put(1,0){\circle*{.3}}
\put(2,0){\circle*{.3}}
\put(4,0){\circle*{.3}}
\put(5,0){\circle*{.3}}
\put(5.6,0.5){\circle*{.3}}
\put(5.6,-0.55){\circle*{.3}}
\put(-.25, -.2){$\times$}

\put(0,1.5){\line(1,0){1}}
\put(1,1.5){\line(1,0){.85}}
\put(2.1,1.5){\line(1,0){.4}}
\multiput(2.65,1.5)(0.2,0){5}{\line(1,0){.1}}
\put(3.55,1.5){\line(1,0){.4}}
\put(4,1.5){\line(1,0){1}}
\put(5,1.5){\line(1,1){0.6}}
\put(5,1.5){\line(1,-1){0.6}}

\put(1,1.5){\circle*{.3}}
\put(2,1.5){\circle*{.3}}
\put(4,1.5){\circle*{.3}}
\put(5,1.5){\circle*{.3}}
\put(5.6,2.1){\circle*{.3}}
\put(5.6,1.){\circle*{.3}}
\put(-.25, 1.3){$\times$}
\put(3, 0.6){$\updownarrow$$\sigma$}
\end{picture} \\

\textrm{DIII}& $\PO_{4n}$&$\PGL_{2n}$ &$\Gr(2,2n)$ & 
\begin{picture}(8,3)(0,0.5)
\put(0,0){\line(1,0){1}}
\put(1,0){\line(1,0){.85}}
\put(2.1,0){\line(1,0){.4}}
\multiput(2.65,0)(0.2,0){5}{\line(1,0){.1}}
\put(3.55,0){\line(1,0){.4}}
\put(4,0){\line(1,0){1}}
\put(5,0){\line(1,0){1}}

\put(0,0){\circle*{.3}}
\put(2,0){\circle*{.3}}
\put(4,0){\circle*{.3}}
\put(5,0){\circle*{.3}}
\put(6,0){\circle*{.3}}
\put(0.7, -0.2){$\times$}

\put(0,1){\line(1,0){1}}
\put(1,1){\line(1,0){.85}}
\put(2.1,1){\line(1,0){.4}}
\multiput(2.65,1)(0.2,0){5}{\line(1,0){.1}}
\put(3.55,1){\line(1,0){.4}}
\put(4,1){\line(1,0){1}}
\put(5,1){\line(1,0){1}}

\put(0,1){\circle*{.3}}
\put(2,1){\circle*{.3}}
\put(4,1){\circle*{.3}}
\put(1,1){\circle*{.3}}
\put(6,1){\circle*{.3}}
\put(4.7, 0.8){$\times$}
\put(6.5, 0.2){$\updownarrow$$\sigma$}
\end{picture}
\\&&& $\sqcup \Gr(2n-2,2n)$ & \\

\textrm{EVII}& $E^{ad}_7$ & $E_6 \times \mathbb C^*$ & $E_6/P_1 \sqcup E_6/P_6$ &

\begin{picture}(8,1.5)(0,1.5)
\put(1,0){\line(1,0){1}}
\put(2,0){\line(1,0){1}}
\put(3,0){\line(1,0){1}}
\put(4,0){\line(1,0){1}}
\put(3,0){\line(0,-1){0.8}}
\put(2,0){\circle*{.3}}
\put(3,0){\circle*{.3}}
\put(4,0){\circle*{.3}}
\put(5,0){\circle*{.3}}
\put(3,-1){\circle*{.3}}
\put(0.75,-.15){$\times$}

\put(1,2){\line(1,0){1}}
\put(2,2){\line(1,0){1}}
\put(3,2){\line(1,0){1}}
\put(4,2){\line(1,0){1}}
\put(3,2){\line(0,-1){0.8}}
\put(2,2){\circle*{.3}}
\put(3,2){\circle*{.3}}
\put(4,2){\circle*{.3}}
\put(1,2){\circle*{.3}}
\put(3,1){\circle*{.3}}
\put(4.75,1.85){$\times$}
\put(5.5, 0.7){$\updownarrow$$\sigma$}
\end{picture} \\
&&&&\\ &&&&\\
&&&&\\ &&&&\\
\hline
    \end{tabular}
    \caption {Hermitian exceptional and non-exceptional type $G/N_G(L)$}
\end{table}

\clearpage

\vspace{1em}
\noindent
{\bf Acknowledgements.} The author appreciate to Michel Brion and Nicolas Perrin for helpful comments and encourage. The author thanks the Institut Fourier for hosting her stay from 2018 to 2019 which led to this work. She also thanks the Center for Geometry and Physics of Institute for Basic Science and Basic Science Research Institute of Ewha Womans University for their support.

\vspace{1em}

\noindent
\textbf{Funding}. 
This work was partially supported by the Institute for Basic Science (IBS-R003-D1) and by Basic Science Research Program through NRF Korea (2021R1A6A1A10039823).

\vspace{1em}

\noindent
\textbf{Data availability}. 
Data sharing not applicable to this article as no datasets were generated or analysed during the current study.

\vspace{1em}

\noindent
\textbf{Competing interests}. 
The authors have no competing interests to declare that are relevant to the content of this article.

\end{document}